\documentclass[12pt]{amsart}

\usepackage{color}

\usepackage{amssymb}
\usepackage{graphicx}  
\DeclareGraphicsExtensions{.pstex,.eps}

\newcommand{\CC}{{\mathbb C}}

\newcommand{\RR}{{\mathbb R}}

\newcommand{\defeq}{\stackrel{\rm{def}}{=}}

\renewcommand{\l}{\lambda}
\newcommand{\supp}{\operatorname{supp}}
\newcommand{\kerr}{\operatorname{ker}}

\newcommand{\tr}{\operatorname{tr}}

\newcommand{\Res}{\operatorname{Res}}
\newcommand{\rest}{\!\!\restriction}

\renewcommand{\Re}{\mathop{\rm Re}\nolimits}
\renewcommand{\Im}{\mathop{\rm Im}\nolimits}

\theoremstyle{plain}

\newtheorem{thm}{Theorem}
\newtheorem{prop}{Proposition}[section]

\newtheorem{lem}[prop]{Lemma}

\theoremstyle{definition}

\numberwithin{equation}{section}

\def\bbbone{{\mathchoice {1\mskip-4mu {\rm{l}}} {1\mskip-4mu {\rm{l}}}
{ 1\mskip-4.5mu {\rm{l}}} { 1\mskip-5mu {\rm{l}}}}}

\def\squarebox#1{\hbox to #1{\hfill\vbox to #1{\vfill}}} 
 
\newcommand{\sech}{\textnormal{sech}}

\newcommand{\indentalign}{\hspace{0.3in}&\hspace{-0.3in}}
\newcommand{\la}{\langle}
\newcommand{\ra}{\rangle}

\newcommand{\nlso}{\textnormal{NLS}_0}
\newcommand{\nlsq}{\textnormal{NLS}_q}
\newcommand{\trans}{\textnormal{tr}}
\newcommand{\refl}{\textnormal{ref}}

\usepackage{amsxtra}

\title
[Fast soliton scattering by delta impurities]
{Fast soliton scattering by delta impurities}

\author[J. Holmer]
{Justin Holmer}
\author[J. Marzuola]
{Jeremy Marzuola}
\author[M. Zworski]
{Maciej Zworski}
\address{Mathematics Department, University of California \\
Evans Hall, Berkeley, CA 94720, USA}

\begin{document}    
   
\maketitle   

\begin{abstract}
We study the Gross-Pitaevskii equation with a repulsive delta function potential.  We show that a high velocity incoming soliton is split into a transmitted component and a reflected component.   The transmitted mass ($L^2$ norm squared) is shown to be in good agreement with the quantum transmission rate of the delta function potential.  We further show that the transmitted and reflected components resolve into solitons plus dispersive radiation, and quantify the mass and phase of these solitons.
\end{abstract}
   
\section{Introduction}   
\label{in}

We study the Gross-Pitaevskii equation (NLS) with a repulsive delta function potential ($q>0$)
\begin{equation}
\label{eq:nls}
\left\{
\begin{aligned}
&i\partial_t u + \tfrac{1}{2}\partial_x^2 u -q\delta_0(x)u +u|u|^2 = 0\\
&u(x,0) = u_0(x)
\end{aligned}
\right.
\end{equation}
As initial data we take a fast soliton approaching the impurity from the left:
\begin{equation}
\label{eq:init}
 u_0 ( x ) = e^{ i v x } \sech ( x - x_0 ) \,, \ \ v \gg 1 \,,
\ \  x_0 \ll 0.
\end{equation}
Because of the homogeneity of the problem this covers the case of the general soliton profile $ A \sech ( A x ) $.  The quantum transmission rate at velocity $ v $ is given by the square of the absolute value of the transmission coefficient, 
see  \eqref{eq:tr} below,
\begin{equation}
\label{eq:tqv}
  T_q ( v ) = | t_q ( v ) |^2 = \frac{ v^2}{ v^2 + q^2 } \,.
\end{equation}
For the soliton scattering the natural definition of the transmission rate is given by 
\begin{equation}
\label{eq:tqs} 
 T_q^{\rm{s}} ( v ) = \lim_{ t\rightarrow \infty } 
\frac{ \| u ( t ) \rest_{ x > 0 } \|_{L^2} ^2 }{ \| u ( t )  \|_{ L^2}^2} 
 =  \frac{1} 2  \lim_{ t\rightarrow \infty } 
{ \| u ( t ) \rest_{ x > 0 } \|_{L^2} ^2 }
\,,\end{equation}
provided that the limit exists.   We expect that it {\em does} and that for fixed $q/v$, there is a $\sigma>0$ such that 
\begin{equation}
\label{eq:conj}
T^{\rm{s}}_q(v) = T_q(v) + {\mathcal O}(v^{-\sigma}),   
\quad \text{as }v\to +\infty  \,.
\end{equation}
Based on the comparison with the linear case (see \eqref{eq:rem} below)
and the numerical evidence \cite{HMZ2} we expect \eqref{eq:conj} with $ \sigma = 2$.

\begin{figure}
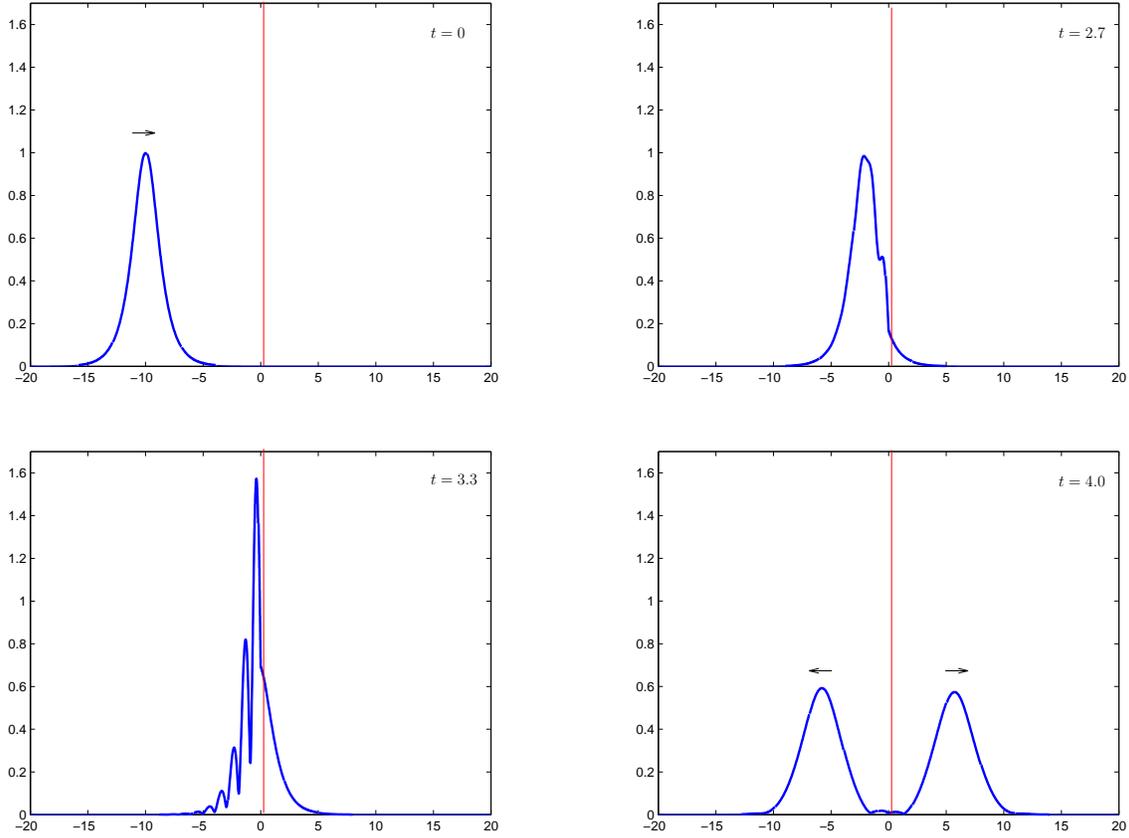

\scalebox{0.5}{\input{snap1a.pstex_t}} \hfill \scalebox{0.5}{\input{snap1b.pstex_t}}
\scalebox{0.5}{\input{snap1c.pstex_t}} \hfill \scalebox{0.5}{\input{snap1d.pstex_t}}

\caption{\label{f:solv3} Numerical simulation of the case $q=v=3$, $x_0=-10$, at times $t=0.0, 2.7, 3.3, 4.0$.  Each frame is a plot of amplitude $|u|$ versus $x$.}
\end{figure}

Towards this heuristic claim we have

\begin{thm}
\label{th:1}
Let $ \delta $ satisfy $\frac{2}{3}<\delta<1$.  
If $u(x,t)$ is the solution of \eqref{eq:nls} with initial condition 
\eqref{eq:init} and $x_0\leq -v^{1-\delta}$, then for fixed $q/v$,
\begin{equation}
\label{eq:th}
 \frac12 \int_{x>0} |u(x,t)|^2 \, dx = \frac{v^2}{v^2+q^2} 
+ \mathcal{O}(v^{1-\frac32\delta}), \quad \text{as }v \to +\infty\,, 
\end{equation}
uniformly for 
$$\frac{|x_0|}{v} + v^{-\delta} \leq t \leq (1-\delta)\log v$$
\end{thm}

We see that by taking $\delta$ very close to $1$, we obtain an asymptotic rate 
just shy of $v^{-1/2}$.  More precisely, we show that there exists 
\[ v_0 = v_0 ( q/v , \delta ) \,,\]
diverging   to $+\infty$ as $\delta \uparrow 1$ and $q/v \to +\infty$, 
such that for fixed $ q/v $, if $v\geq v_0$, then
$$\left| \frac{1}{2} \int_{x>0} |u(x,t)|^2 dx - \frac{v^2}{v^2+q^2} \right| \leq cv^{1-\frac{3}{2}\delta} .$$
The constant $c$ appearing here is 
independent of all parameters ($q$, $v$, and $\delta$).  

\begin{figure}
\scalebox{0.8}{\input{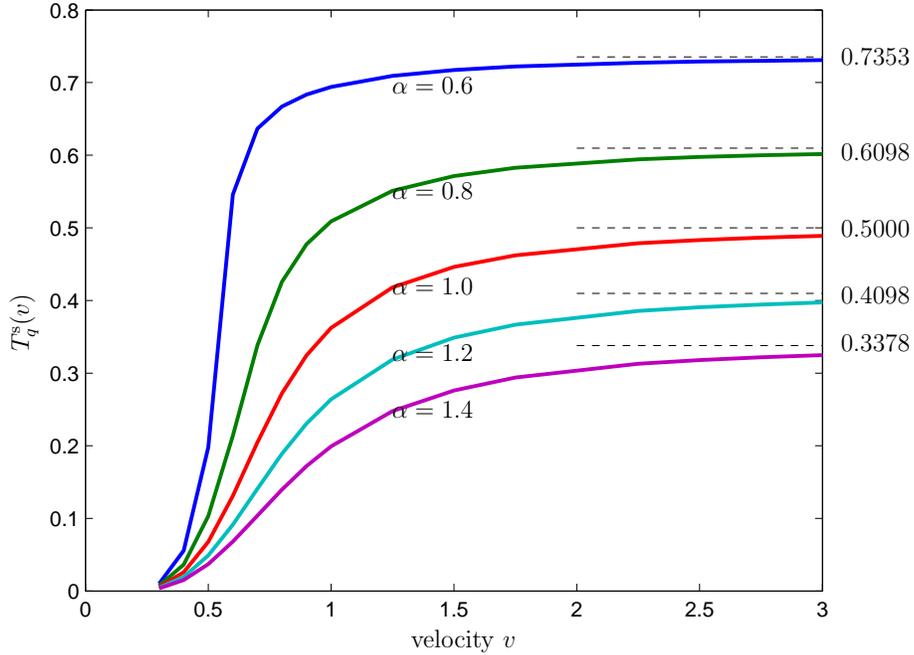}}
\caption{\label{f:nice} A plot of the numerically obtained transmission $T_q^{\text{s}}(v)$ versus velocity $v$ for five values of $\alpha =q/v=0.6,0.8,1.0,1.2,1.4$.  The dashed lines are the corresponding theoretical $v\to +\infty$ asymptotic values given by $1/(1+\alpha^2)$.}
\end{figure}

We have conducted a numerical verification of Theorem \ref{th:1} -- see Fig.\ \ref{f:nice}.  It shows that the approximation given by \eqref{eq:th} is very good even for velocities as low as $\sim 3$, at least for 
\[ 0.6 \leq \alpha \defeq q/v \leq 1.4 \,.\] 
A more elaborate numerical analysis will appear 
in our forthcoming paper \cite{HMZ2}.

Our second result shows that the scattered solution is given, on the same time scale, by a sum of a reflected and a transmitted soliton, and of a time decaying (radiating) term -- see the fourth frame of Fig.\ \ref{f:solv3}. This is further supported by a forthcoming numerical study \cite{HMZ2}.  In previous works in the physics literature (see for instance \cite{CM}) the resulting waves were only described as ``soliton-like''.

\begin{thm}
\label{th:2} 
Under the hypothesis of Theorem \ref{th:1} and for 
$$ \frac{|x_0|}{v}+1 \leq  t  \leq (1-\delta) \log v, $$ 
we have, as $v\to +\infty$, 
\begin{gather}
\label{eq:th2}
\begin{gathered}
u(x,t)  = u_T ( x , t) + u_R ( x , t ) 
 + \mathcal{O}_{L_x^\infty}\left(\left(t-{|x_0|}/{v}\right)^{-1/2}\right) 
+ {\mathcal O}_{L_x^2}( v^{1-\frac{3}{2}\delta})  \,,  \\
 u_T ( x , t )  =  e^{ i \varphi_T} 
e^{ixv + i  ( A_T^2 - v^2 )t /2 }  A_T \, \sech(A_T(x-x_0-tv)) \,, \\
u_R ( x, t)   =  e^{i \varphi_R} 
 e^{-ixv + i  ( A_R^2 - v^2 ) t/2 } A_R \, \sech(A_R(x+x_0+tv)) \,, 
\end{gathered}
\end{gather}
where $ A_T = (2|t_q(v)|-1)_+  $, $ A_R = (2|r_q(v)|-1)_+ $, and
\begin{gather}
\label{eq:th3}
\begin{gathered}
\varphi_T = \arg t_q ( v ) + \varphi_0(|t_q(v)|) + (1-A_T^2)|x_0|/2v \,, 
\\
\varphi_R = \arg r_q ( v ) + \varphi_0(|r_q(v)|) + (1-A_R^2)|x_0|/2v \,, 
\end{gathered}
\end{gather}
$$\varphi_0(\omega) =  \int_0^\infty \log\left( 1 + \frac{\sin^2\pi \omega}{\cosh^2\pi \zeta} \right) \frac{\zeta}{\zeta^2+(2\omega-1)^2} \, d\zeta \,. $$
Here $ t_q ( v ) $ and $ r_q ( v ) $ are the transmission and reflection coefficients of the delta-potential (see \eqref{eq:tr}).  When $ 2 | t_q ( v ) | =1 $ or $ 2 | r_q ( v ) | = 1 $ the first error term in \eqref{eq:th2} is modified to $ {\mathcal O}_{ L^\infty_x } ( (\log  ( t - |x_0|/v ))/( t - |x_0|/v ))^{\frac12} ) $.  
\end{thm}
Here and later we use the standard notation
\begin{equation}
\label{eq:sta}  a_+^k =
\left\{ \begin{array}{ll} a^k & a \geq 0 \,, \\
  0 & a < 0 \,. \end{array} \right. 
\end{equation}

This asymptotic description holds for $v$ greater than some threshold depending on $q/v$ and $\delta$, as in Theorem \ref{th:1}.  The implicit constant in the $\mathcal{O}_{L_x^2}$ error term is entirely independent of all parameters ($q$, $v$, and $\delta$), although the implicit constant in the $\mathcal{O}_{L_x^\infty}$ error term depends upon $q/v$, or more precisely, the proximity of $|t_q(v)|$ and $|r_q(v)|$ to $\frac{1}{2}$.

\begin{figure}
\begin{center}
\scalebox{0.7}{\input{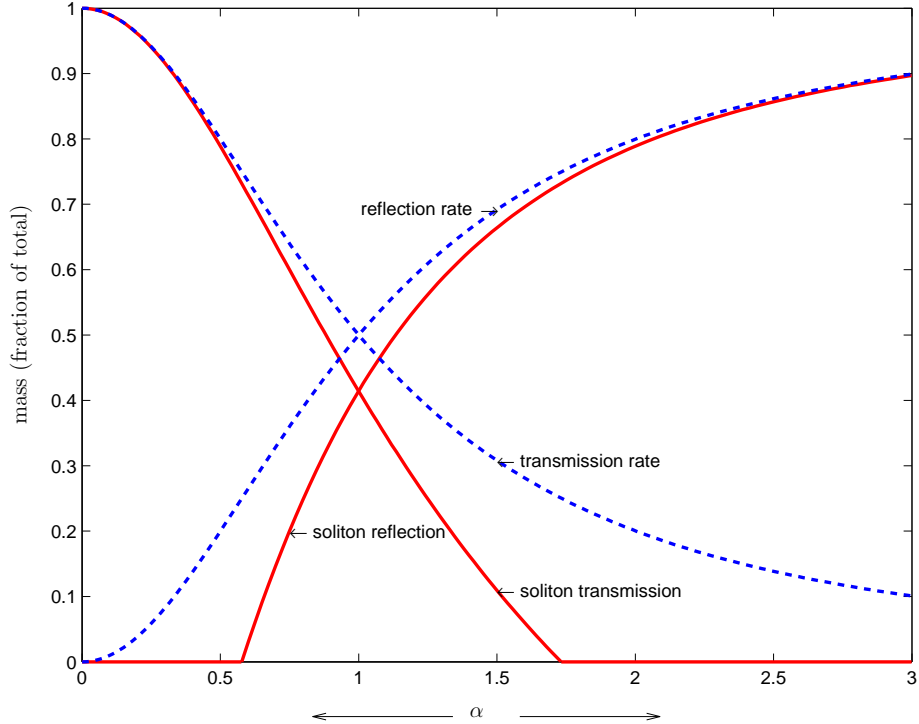}}
\end{center}
\caption
{\label{f:lnl} Comparison of linear and nonlinear scattering 
coefficients as functions of $ \alpha \defeq q / v $.}
\end{figure}

A comparison of the transmission and reflection coefficients \eqref{eq:tqv} of the $ \delta $ potential, and of the soliton transmission and reflections coefficients \eqref{eq:th3}, is shown in Figure \ref{f:lnl}.

Scattering of solitons by delta impurities is a natural model explored extensively in the physics literature -- see for instance \cite{CM},\cite{GHW}, and references given there.  The heuristic insight that at high velocities ``linear scattering'' by the external potential should dominate the partition of mass is certainly present there.  In the mathematical literature the dynamics of solitons in the presence of external potentials has been studied in high velocity or semiclassical limits following the work of Floer and Weinstein \cite{FlWe}, and Bronski and Jerrard \cite{BJ} -- see \cite{FrSi} for recent results and a review of the subject. Roughly speaking, the soliton evolves according to the classical motion of a particle in the external potential. That is similar to the phenomena in other settings, such as the motion of the Landau-Ginzburg vortices.

The possible novelty in \eqref{eq:th} and \eqref{eq:th2} lies in seeing {\em quantum} effects of the external potential strongly affecting soliton dynamics. As shown in Fig.\ \ref{f:nice}, Theorem \ref{th:1} gives a very good approximation to the transmission rate already at  low velocities.  Fig.\ \ref{f:solv3} shows time snapshots of the evolution of the soliton, and the last frame suggests the 
soliton resolution \eqref{eq:th2}. We should stress that the asymptotic 
solitons are resolved at a much larger time -- see \cite{HMZ2}.

The proof of the two theorems, given below in \S\ref{proof}--\ref{resolution}, proceeds by approximating the solution during the ``interaction phase'' (the interval of time during which the solution significantly interacts with the delta potential at the origin) by the corresponding linear flow.  This approximation is achieved, uniformly in $q$, by means of Strichartz estimates established in \S\ref{ros}.  The use of the Strichartz estimates as an approximation device, as opposed to say energy estimates, is critical since the estimates obtained depend only upon the $L^2$ norm of the solution, which is conserved and \textit{independent} of $v$.  Thus, $v$ functions as an asymptotic parameter; larger $v$ means a shorter interaction phase and a better approximation of the solution by the linear flow.  Theorem 2 combines this analysis with the inverse scattering method. The delta potential splits the incoming soliton into two waves which become single {\em solitons}.

\medskip

\noindent
{\sc Acknowledgments.} We would like to thank Mike Christ, Percy Deift, 
and Michael Weinstein for helpful discussions during the preparation of this paper. The work of the first author was supported in part by an NSF postdoctoral fellowship, and that
of the second and third author by NSF grants  DMS-0354539 and DMS-0200732.

\section{Scattering by a delta function}
\label{ros}

Here we present some basic facts about scattering by a delta-function potential on the real line.  Let $ q \geq 0  $ and put
\[ H_q = - \frac{1}2 \frac{d^2}{dx^2} + q \; \delta_0 ( x ) \,.\]
We define special solutions, $ e_\pm ( x , \lambda ) $,
to $ ( H_q - \lambda^2 /2 ) e_\pm  = 0 $, using notation given in \eqref{eq:sta}:
\begin{equation}
\label{211}
e_{\pm}(x,\l) = t_q (\l)e^{\pm i \l x} x_{\pm}^0 
+ (e^{\pm i \l x} + r_q (\l)e^{\mp i\l x}) x_{\mp}^0  \,, 
\end{equation}
where $ t_q $ and $ r_q $ are the 
the transmission and reflection coefficients:
\begin{equation}
\label{eq:tr}
t_q ( \lambda ) = \frac{ i \lambda } { i \lambda - q } \,, \ \ 
r_q ( \lambda ) = \frac{ q} {i \lambda - q } \,.
\end{equation}
They satisfy two equations, one standard (unitarity) 
and one due to the special structure of the potential:
\begin{equation}
\label{eq:trpr} | t_q ( \lambda ) |^2 + | r_q ( \l ) |^2 = 1 \,, \ \ 
t_q ( \lambda ) = 1 + r_q ( \lambda ) \,.\end{equation}
We use the representation of the propagator in terms of the generalized eigenfunctions-- see for instance the notes \cite{TZ} covering scattering by compactly supported
 potentials.   The resolvent 
\[ R_q ( \lambda ) \defeq ( H_q - \lambda^2 / 2 )^{-1} \,,\]
has kernel given by 
\[R_q ( \l)(x,y) =\frac{1}{i\l t_q (\l)}\,
\big(e_+(x,\l)e_-(y,\l)(x-y)^0_+ +
  e_+(y,\l)e_-(x,\l)(x-y)^0_-\big) \,.\]
This gives an explicit formula for the spectral projection, and hence the Schwartz kernel of the propagator:
\begin{equation}
\label{eq:sppr}
\exp ( - i t H_q )  = \frac{1}{2\pi}\int^\infty_0 
e^{- i t \lambda^2/2 } 
 \left(e_+(x,\l)\overline{e_+(y,\l)} + e_-(x,\l)  
\overline{e_-(y,\l)}\right) \,d\l \,.
\end{equation}
The propagator for $ H_q  $ is described in the following

\begin{lem}
\label{L:lin}
Suppose that $ \varphi \in L^1 $ and that $ \supp \varphi \subset ( -\infty , 0] $.
Then 
\begin{equation}
\label{eq:prop}
\begin{split}
&  \exp ( - i t H_q ) \varphi ( x ) = 
\\ 
& \ \ \ \exp ( - i t H_0 ) ( \varphi * \tau_q ) ( x ) 
x_+^0 + 
( \exp ( - i t H_0 ) \varphi ( x ) + \exp ( - it H_0 ) ( \varphi * \rho_q ) ( - x) 
) x_-^0  \,, 
\end{split}
\end{equation}
where 
\begin{equation}
\label{eq:ftr}  \rho_q ( x) = - q \exp(qx) x_-^0  \,, \ \ \tau_q 
( x  ) = \delta_0 ( x) + \rho_q ( x ) \,. 
\end{equation}
\end{lem}
\begin{proof}
All we need to do is to combine \eqref{211} and \eqref{eq:sppr}.
Using the support property of $ \varphi $ we compute,
\[ \begin{split}
&  \int \varphi ( y ) \overline{e_+ ( y , \lambda ) } dy = r_q ( - \lambda ) 
\hat \varphi ( - \lambda ) + \hat \varphi ( \lambda ) \,, \\
& \int \varphi ( y ) \overline{ e_- ( y , \lambda ) } dy = t_q ( - \lambda ) 
\hat \varphi ( - \lambda )  \,,
\end{split} \]
so that 
\[\begin{split} \left( \exp ( - i t H_q ) \varphi \right) \rest_{ x > 0 } & = 
\frac{1}{2 \pi } \int_0^\infty e^{ - i t \lambda^2/2 } 
\left( t_q ( \lambda ) e^{ i \lambda x } ( r_q ( - \lambda ) 
\hat \varphi ( - \lambda ) + \hat \varphi ( \lambda ) ) \; + \ \right. \\
& \ \ \ \ \ \ \ \ \ \ \ \ \ \ \  
\left. ( r_q ( \lambda ) e^{ i \lambda x } + e^{- i \lambda x } ) t_q ( - 
\lambda ) \hat \varphi ( - \lambda ) \right) d \lambda \\ 
& = \; \frac{1}{ 2 \pi } \int_\RR e^{ - i t \lambda^2/2 } t_q ( \lambda ) 
\hat \varphi ( \lambda ) e^{ i \lambda x } d \lambda \\
& = \exp ( - i t H_0 ) ( \tau_q * \varphi ) ( x ) 
\,, \ \ \widehat \tau_q ( \lambda ) = t_q ( \lambda ) \,,
\end{split} \]
where we used the fact that $ r_q ( - \lambda ) t_q ( \lambda ) + 
r_q ( \lambda ) t_q ( - \lambda ) = 0 $.

Similarly, using $ r_q ( - \lambda ) r_q ( \lambda ) + 
t_q (- \lambda ) t_q ( \lambda ) = 1 $, we have 
\[ \begin{split} \left( \exp ( - i t H_q ) \varphi \right) \rest_{ x < 0 } & = 
\frac{1}{ 2 \pi } \int_0^\infty e^{ - i
 t \lambda^2/2 } \left( 
\hat \varphi ( \lambda ) e^{ i \lambda x } + r_q ( \lambda ) \hat \varphi ( \lambda ) 
e^{ - i \lambda x } \right) d\lambda \\
& = \; \exp( - i t H_0 ) \varphi ( x ) + \exp ( - i t H_0 ) ( \varphi * \rho_q ) ( - x ) 
\,, \ \ \widehat {\rho_q} ( \lambda ) = r_q ( \lambda ) \,. 
\end{split} \] 
A simple computation gives \eqref{eq:ftr} concluding the proof.
\end{proof}

We have two simple applications of Lemma \ref{L:lin}:  the Strichartz estimates  (Proposition \ref{p:Str}) and the asymptotics of the linear flow $\exp(-itH_q)$ as $ v \rightarrow +\infty $ (Proposition \ref{p:as}).  We start with the Strichartz estimate, which will be used several times in the various approximation arguments of \S\ref{proof}.  Since it is particularly simple in our setting, we give a complete proof (see \cite{KT} for references and the general version).
\begin{prop}
\label{p:Str}
Suppose  $q\geq 0$ and 
\begin{equation}
\label{eq:eq}
 i \partial_t u ( x , t ) + \tfrac{1}{2}\partial_{x}^2 u ( x, t ) 
- q \delta_0 ( x ) u ( x , t ) = f ( x , t ) \,, \ \ u ( x , 0 ) = \varphi ( x ) 
\,.\end{equation}
Let the indices $p,r$, $\tilde p$, $\tilde r$ satisfy
\begin{equation}
\label{eq:adm}
2 \leq p, r \leq \infty \,, \ \ 1 \leq \tilde p , \tilde r \leq 2 \,, \ \
\frac 2 p + \frac 1 r = \frac 12 \,, \ \ \ \frac 2 {\tilde p }
+ \frac 1 {\tilde r} = \frac 52
\end{equation}
and fix a time $T>0$.  Then
\begin{equation}
\label{eq:Str}
\| u \|_{ L^p_{[0,T]} L^r_x } \leq c \| \varphi \|_{L^2} + c \| f \|_{ L_{[0,T]}^{\tilde p} L_x^{\tilde r} }
\end{equation}
The constant $c$ is independent of $q$ and $T$.  Moreover, in \eqref{eq:eq}, we can take $f(x,t) = g(t)\delta_0(x)$ and, on the right-hand side of \eqref{eq:Str}, replace $\| f \|_{ L_{[0,T]}^{\tilde p} L_x^{\tilde r} }$  with $\|g\|_{L_{[0,T]}^\frac{4}{3}}$.
\end{prop}
\begin{proof}
We put $ U_q ( t ) \defeq \exp ( - i t H_q ) $, so that $ U_q ( t ) $ is 
a unitary group on $ L^2 (\RR ) $. For $ \varphi \in L^1 ( \RR ) $ we have,
using Lemma \ref{L:lin},
\begin{equation}
\label{eq:utl}
\begin{split} 
\| U_q ( t ) \varphi \|_{ L^\infty } & \leq \sum_{\pm } \| U_q ( t ) ( \varphi x_\pm^0 ) 
\|_{L^\infty } \\
& \leq  \sum_{ \pm } \| U_0 ( t ) \|_{ L^1 \rightarrow L^\infty } ( \| 
( \varphi x_\pm^0 ) * \tau_q \|_{L^1 } + 
 \| ( \varphi x_\pm^0 ) * \rho_q \|_{L^1 } ) \\ 
& \leq \frac{ 1}{ \sqrt {\pi |t|} }  ( 1 + 2\| \rho_q \|_{L^1} )
 \| \varphi\|_{L^1} \\
& \leq \frac{3}{\sqrt{ \pi |t|}} \| \varphi \|_{ L^1} \,. 
\end{split}
\end{equation}
By the Riesz-Thorin interpolation theorem (see for instance \cite[Theorem 7.1.12]{Hor1}) we have
\begin{equation}
\label{eq:rith} \| U ( t )  \|_{ L^{r'} \rightarrow L^r } 
\leq C | t |^{ -\frac 1 2 \left( 1 - \frac 2 r \right)   } \,, \ \ 
1 \leq r' \leq 2 \,, \ \ \frac 1 r + \frac 1 {r'}= 1\,. 
\end{equation}

The estimate \eqref{eq:Str} with $ f \equiv 0 $ reads
\[ \| U ( t ) g  \|_{ L^p_t L^r_x } \leq 
C \| g \|_{L^2 ( \RR)}  \,, \]
which by duality is equivalent to 
\begin{equation}
\label{eq:Str1}
\left\| \int_\RR U ( - s ) F ( s) ds \right\|_{L^2 ( \RR ) } 
\leq C \| F \|_{ L^{p'}_t L^{r'}_x } \,.
\end{equation}
The two equivalent estimates together give (\eqref{eq:Str1} is applied with $p'$, $r'$ replaced by $\tilde p$, $\tilde r$ -- it is easily checked that \eqref{eq:adm} still holds)
\[  \left\| \int_\RR U ( t - s ) F ( s) ds \right\|_{L^p_t L^r_x }
 \leq  \| F \|_{ L^{\tilde p}_t L^{\tilde r}_x } \,,\]
Putting
$ F ( s ) = \bbbone_{ [0 ,t] } ( s )  f ( s, x ) $ we obtain  \eqref{eq:Str} for $ u_0 = 0 $. Hence it suffices to prove \eqref{eq:Str1}. 

Put 
\[ T F ( x ) \defeq \int_\RR U ( - s ) F ( s , x ) d s \,,\]
and note that $ T^*  g (s,  x ) := U ( s) g ( x ) $.
The estimate \eqref{eq:Str1} is equivalent to
\[  \langle T^* T G , F \rangle_{ L^2_{t,x} } \leq C
 \| G \|_{  L^{p'}_t L^{r'}_x } \| F \|_{  L^{p'}_t L^{r'}_x   } \,,\]
which is the same as
\begin{equation}
\label{eq:utus0}
\left| \int_\RR \! \! \!  \int_\RR  \langle G ( t ) ,\
 U ( t - s ) F ( s ) \rangle
\; dt ds \right| 
 \leq C \| G \|_{ L^{p'}_t L^{r'}_x  }
\| F \|_{  L^{p'}_t L^{r'}_x   } \,.
\end{equation}

To obtain \eqref{eq:utus0} from 
\eqref{eq:rith} we apply 
the Hardy-Littlewood-Sobolev inequality 
which says that if $ K_a ( t ) = | t|^{-1/a } $ and $ 1 < a < \infty  $ 
then 
\begin{gather*}
 \| K_a * F \|_{L^{\alpha}  ( \RR ) } \leq C \| F \|_{ L^{\beta } ( \RR ) } 
\,, \ \ \frac1\alpha = \frac1\beta - \frac1a \,, \ \ 
1 < \beta < \alpha  \,, 
\end{gather*} 
see for instance \cite[Theorem 4.5.3]{Hor1}. We apply it with 
\[  \frac 1 a = \frac12 \left( 1 - \frac 2 r  \right) \,,  \quad  \alpha =p \,, \quad \beta = p' \,, \]
which is the admissibility condition \eqref{eq:adm}.
\end{proof}

We now turn to the large velocity asymptotics of the linear flow $\exp(-itH_q)$.\begin{prop}
\label{p:as}
Let $\theta \in C^\infty ( \RR ) $  be bounded , together will  all of 
its derivatives.  Let $\varphi\in \mathcal{S}(\mathbb{R})$, $v>0$, and suppose 
$\supp  [\theta( \bullet )\varphi( \bullet -x_0)] 
\subset (-\infty,0]$.  Then for $2|x_0|/v \leq t \leq 1$,
\begin{equation}
\label{E:as2}
e^{-itH_q}[e^{ixv}\varphi(x-x_0)] =
\begin{aligned}[t]
& t(v)e^{-itH_0}[e^{ixv}\varphi(x-x_0)] \\
&+ r(v) e^{-itH_0}[e^{-ixv}\varphi(-x-x_0)] \\
&+ e(x,t)
\end{aligned}
\end{equation}
where, for any $k\geq 0$,
$$\|e(\cdot,t)\|_{L^2} \leq 
\begin{aligned}[t]
&\frac{1}{v}\|\partial_x [\theta(x)\varphi(x-x_0)]\|_{L^2}\\
&+\frac{c_k}{(tv)^k} \| \la x \ra^k \varphi(x) \|_{H^k}\\
&+4\|(1-\theta(x))\varphi(x-x_0)\|_{L_x^2}
\end{aligned}$$
\end{prop}
In \S \ref{proof}, Proposition \ref{p:as} will be applied with $\theta(x)$ a smooth cutoff to $x<0$, and $\varphi(x)=\sech x$ with $x_0=-v^{1-\delta} \ll 0$.  

Before proving Proposition \ref{p:as}, we need the following
\begin{lem}
\label{p:as0}
Let $\psi\in \mathcal{S}(\mathbb{R})$ with $\supp \psi \subset (-\infty,0]$.  Then
\begin{equation}
\label{E:as1}
e^{-itH_q}[e^{ixv}\psi(x)](x) =
\begin{aligned}[t]
&e^{-itH_0}[e^{ixv}\psi(x)](x)\, x_-^0 \\
&+t(v)e^{-itH_0}[e^{ixv}\psi(x)](x) \, x_+^0 \\
&+r(v)e^{-itH_0}[e^{-ixv}\psi(-x)](x) \, x_-^0\\
&+e(x,t)
\end{aligned}
\end{equation}
where 
$$\|e(x,t)\|_{L_x^2} \leq \frac{1}{v}\|\partial_x \psi\|_{L^2}$$
uniformly in $t$.
\end{lem}

\begin{proof}[Proof of Lemma \ref{p:as0}]
By \eqref{eq:prop} with $\varphi(x)=e^{ixv}\psi(x)$,
$$e(x,t) =
\begin{aligned}[t]
&[e^{-itH_0}(\varphi\ast (\tau - t(v)\delta_0))(x)]\, x_+^0 \\
&+[e^{-itH_0}(\varphi\ast (\rho - r(v)\delta_0))(-x)]\, x_-^0 
\end{aligned}
$$
and thus it suffices to show
\begin{equation}
\label{E:101}
\|e^{-itH_0}(\varphi\ast (\tau - t(v)\delta_0))(x)\|_{L_x^2} \leq \frac{1}{v}\|\partial_x \psi\|_{L_x^2}
\end{equation}
and
$$\|e^{-itH_0}(\varphi\ast (\rho - r(v)\delta_0))(x)\|_{L_x^2} \leq \frac{1}{v}\|\partial_x \psi\|_{L_x^2}$$

The proofs of these two estimates are similar, so we only carry out the proof of \eqref{E:101}. By unitarity of $ e^{ - i t H _0 } $ and 
the Plancherel's identity, 
\begin{equation}
\label{E:100}
\|e^{itH_0}[\varphi\ast (\tau - t(v)\delta_0)](x)\|_{L_x^2} = 
\| \hat{\psi}(\lambda -v)( t(\lambda)-t(v))\|_{L_\lambda^2} \,. 
\end{equation}
Since
$$ t (\lambda)-t(v) = \frac{-iq(\lambda-v)}{(i\lambda-q)(iv-q)}$$
we have $| t(\lambda)-t(v)| \leq (\lambda -v )/v$.  Using this to estimate the right-hand side of \eqref{E:100} and applying Plancherel's 
identity again yields \eqref{E:101}.
\end{proof}

\begin{proof}[Proof of Proposition \ref{p:as}]
Apply \eqref{E:as1} to $\psi(x)=\theta(x)\varphi(x-x_0)$ to obtain
\begin{equation}
\label{E:102}
e^{-itH_q}[e^{ixv}\varphi(x-x_0)](x) =
\begin{aligned}[t]
&e^{-itH_0}[e^{ixv}\varphi(x-x_0)](x)\, x_-^0 \\
&+t(v)e^{-itH_0}[e^{ixv}\varphi(x-x_0)](x) \, x_+^0 \\
&+r(v)e^{-itH_0}[e^{-ixv}\varphi(-x-x_0)](x) \, x_-^0\\
&+e_1(x,t)+e_2(x,t)
\end{aligned}
\end{equation}
where $e_1(x,t)$ is as in Lemma \ref{p:as0} and (putting $f(x)=e^{ixv}(1-\theta(x))\varphi(x-x_0)$)
$$e_2(x,t) = 
\begin{aligned}[t]
&+e^{-itH_q}f(x) - e^{-itH_0}f(x) \, x_-^0 \\
&-t(v)e^{-itH_0}f(x)\, x_+^0 - r(v) e^{-itH_0}[f(-x)](x) \, x_-^0
\end{aligned}
$$
By Lemma \ref{p:as0},
$$\|e_1(x,t)\|_{L_x^2} \leq \frac{1}{v}\|\partial_x[\theta(x)\varphi(x-x_0)]\|_{L_x^2}$$
uniformly for all $t$, and by unitarity of the linear flows,
$$\|e_2(x,t) \|_{L_x^2} \leq 4\|(1-\theta(x))\varphi(x-x_0)\|_{L_x^2}$$
also uniformly in all $t$.  Now restrict to the time interval $2|x_0|/v \leq t \leq 1$.  By \eqref{E:102}, it remains to show that 
\begin{gather}
\label{E:103} \|e^{-itH_0}[e^{ixv}\varphi(x-x_0)](x)\|_{L_{x<0}^2} \leq \frac{c_k}{(tv)^k}\|\la x \ra^k\varphi(x)\|_{H_x^k} \\
\notag\|e^{-itH_0}[e^{-ixv}\varphi(-x-x_0)](x)\|_{L_{x>0}^2} \leq \frac{c_k}{(tv)^k}\|\la x \ra^k\varphi\|_{H^k}
\end{gather}
The second of these is in fact equivalent to the first, since for any function $g(x)$, 
$$e^{-itH_0}[g(-x)](x) = e^{-itH_0}[g(x)](-x)\,.$$  
Now we establish \eqref{E:103}.  Since 
$${ [{e^{i \bullet v}\varphi( \bullet -x_0) } ]} \, \hat{} \,  (\lambda) 
= e^{-ix_0(\lambda-v)}\hat{\varphi}(\lambda-v) \,, $$  
\begin{align*}
 e^{-itH_0}[e^{ixv}\varphi(x-x_0)](x) 
&=\frac{1}{2\pi} \int e^{ix\lambda} e^{-ix_0(\lambda-v)}e^{-it\lambda^2/2}\hat{\varphi}(\lambda-v)\, d\lambda\\
&=e^{-itv^2/2}e^{ixv}\frac{1}{2\pi} \int e^{i\lambda(x-x_0-tv)}e^{-it\lambda^2/2}\hat{\varphi}(\lambda) \, d\lambda
\end{align*}
By $k$ applications of integration by parts in $\lambda$,
$$\int e^{i\lambda(x-x_0-tv)}e^{-it\lambda^2/2}\hat{\varphi}(\lambda) \, d\lambda = \left( \frac{i}{x-x_0-tv} \right)^k \int e^{i\lambda(x-x_0-tv)}\partial_\lambda^k[e^{-it\lambda^2/2}\hat{\varphi}(\lambda)] \, d\lambda$$
Since $2|x_0|/v \leq t$, we have $-x_0-tv<0$ and thus $|x-x_0-tv| \geq |-x_0-tv| \geq tv/2$ for $x<0$.  Hence
\begin{align}
\label{eq:hence}
\indentalign \int_{-\infty}^0 |e^{-itH_0} [e^{ixv}\varphi(x-x_0)](x)|^2 \, dx 
\nonumber \\
&\leq \frac{c_k}{(tv)^k} \left\| \int e^{i\lambda(x-x_0-tv)}\partial_\lambda^k[e^{-it\lambda^2/2}\hat{\varphi}(\lambda)] \, d\lambda \right\|_{L_x^2} \\
& =  \frac{c_k}{(tv)^k} \left\| \partial_\lambda^k[e^{-it\lambda^2/2}\hat{\varphi}(\lambda)] \right\|_{L_\lambda^2} \nonumber 
\end{align}
from which the result follows by applying the Leibniz product rule and 
the Plancherel identity once again (and using that $t\leq 1$).
\end{proof}

\medskip

\noindent{\bf Remark.} Suppose that 
\[  u ( x , t ) = e^{ - i t H_q } [e^{ixv}\psi ( x ) ] \,, \ \ 
\psi \in {\mathcal S} ( \RR ) \,, \ \ \supp \psi \subset ( - \infty , 0) 
 \,, \ \ \| \psi \|_{L^2 } = 1 \,.\]
Then for $ t \gg 1 $ and as $ v \rightarrow + \infty $,
\begin{equation}
\label{eq:rem} \int_0^\infty | u ( x , t ) |^2 dx = \frac{v^2}{ v^2 + q^2} 
+ {\mathcal O} \left( \frac 1 {v^2} \right) \,. 
 \end{equation}
In fact using \eqref{eq:prop} and an estimate similar to \eqref{eq:hence}
we see that for $ t \geq 1 $
\[ \begin{split} 
\int_0^\infty | u ( x , t ) |^2 dx & = \| e^{ - i t H_0 } ( ( e^{ i \bullet 
v } \psi ) * \tau_q  ) x_+^0 \|^2 = \|  ( e^{ i \bullet 
v } \psi ) * \tau_q \|_2^2 + {\mathcal O } ( v^{-\infty } ) \\
& = \frac{1}{ 2 \pi } \| i \lambda \hat \psi ( \lambda - v ) / ( 
i \lambda - q ) \|^2 + {\mathcal O } ( v^{-\infty } ) \\
& = \frac{1}{ 2 \pi } \int_{ |\lambda - v | \leq \sqrt v } 
 \frac{  \lambda^2 } { \lambda^2 + q^2 } 
| \hat \psi ( \lambda - v ) |^2 d \lambda  + {\mathcal O } ( v^{-\infty } )  
\end{split} \] 
An expansion in powers of $ ( \lambda - v )  /v $ gives \eqref{eq:rem}.

\section{Soliton scattering}
\label{proof}

In this section, we prove Theorem \ref{th:1}.  We recall the notation for operators from Sect.\ref{ros} and
introduce short hand notation for the nonlinear flows:
\begin{itemize}
\item $H_0=-\frac{1}{2} \partial_x^2$.  The flow $e^{-itH_0}$ is termed the ``free linear flow''
\medskip
\item $H_q = -\frac{1}{2} \partial_x^2+q\delta_0(x)$.  The flow $e^{-itH_q}$ is termed the ``perturbed linear flow''
\medskip
\item $\nlsq(t)\varphi$, termed the ``perturbed nonlinear flow'' is the evolution of initial data $\varphi(x)$ according to the equation $i\partial_tu + \tfrac{1}{2}\partial_x^2 u - q\delta_0(x)u + |u|^2u=0$
\medskip
\item $\nlso(t)\varphi$, termed the ``free nonlinear flow'' is the evolution of initial data $\varphi(x)$ according to the equation $i\partial_th + \tfrac{1}{2}\partial_x^2 h + |h|^2h=0$
\end{itemize}

From Sect.\ref{in} we recall the form of the initial condition:  $u_0(x) = e^{ixv}\sech(x-x_0)$, $v \gg 1$, $x_0 \leq -v^{1-\delta}$, $\frac{2}{3}<\delta<1$, and we put $u(x,t) = \nlsq(t)u_0(x)$.

We begin by outlining the scheme, and will then supply the details.  The $\mathcal O$ notation always means $L_x^2$ difference, uniformly on the time interval specified, and up to a multiplicative factor that is independent of $q$, $v$, and $\delta$ (any such dependence will be exhibited explicitly).

\noindent\textbf{Phase 1 (Pre-interaction)}.  Consider $0\leq t \leq t_1$, where $t_1 = |x_0|/v - v^{-\delta}$ so that $x_0+vt_1=-v^{1-\delta}$.  The soliton has not yet encountered the delta obstacle and propagates according to the free nonlinear flow
\begin{equation}
\label{E:approx1}
u(x,t) = e^{-itv^2/2}e^{it/2}e^{ixv}\sech(x-x_0-vt) + \mathcal{O}(qe^{-v^{1-\delta}}), \quad 0\leq t\leq t_1
\end{equation}
The analysis here is valid provided $v$ is greater than some absolute threshold (independent of $q$, $v$, or $\delta$).  But if we further require that $v$ be sufficiently large so that $v^{-3/2}e^{v^{1-\delta}} \geq \alpha$ (recall $\alpha=q/v$), then $qe^{-v^{1-\delta}} \leq v^{-1/2}\leq v^{-\delta/2}$.  This is the error that arises in the main argument of Phase 2 below.

\noindent\textbf{Phase 2 (Interaction)}.  Let $t_2 = t_1+2v^{-\delta}$ and consider $t_1\leq t \leq t_2$.  The incident soliton, beginning at position $-v^{1-\delta}$, encounters the delta obstacle and splits into a transmitted component and a reflected component, which by time $t=t_2$, are concentrated at positions $v^{1-\delta}$ and $-v^{1-\delta}$, respectively.  More precisely, at the conclusion of this phase (at $t=t_2$),
\begin{equation}
\label{E:approx4}
u(x,t_2) = 
\begin{aligned}[t]
&t(v)e^{-it_2v^2/2}e^{it_2/2}e^{ixv}\sech(x-x_0-vt_2)\\
&+r(v)e^{-it_2v^2/2}e^{it_2/2}e^{-ixv}\sech(x+x_0+vt_2) \\
&+ \mathcal{O}(v^{-\frac{1}{2}\delta})
\end{aligned}
\end{equation}

This is the most interesting phase of the argument, which proceeds by using the following three observations
\begin{itemize}
\item The perturbed nonlinear flow is approximated by the perturbed linear flow for $t_1\leq t \leq t_2$.
\item The perturbed linear flow is split as the sum of a transmitted component and a reflected component, each expressed in terms of the free linear flow of soliton-like waveforms.
\item The free linear flow is approximated by the free nonlinear flow on $t_1\leq t \leq t_2$.  Thus, the soliton-like form of the transmitted and reflected components obtained above is preserved.
\end{itemize}
The brevity of the time interval $[t_1,t_2]$ is critical to the argument, and validates the approximation of linear flows by nonlinear flows.

\noindent\textbf{Phase 3 (Post-interaction)}.  Let $t_3=t_2+(1-\delta)\log v$, and consider $[t_2,t_3]$.  The transmitted and reflected waves essentially do not encounter the delta potential and propagate according to the free nonlinear flow, 
\begin{equation}
\label{E:post}
u(x,t) = 
\begin{aligned}[t]
& e^{-itv^2/2}e^{it_2/2}e^{ixv}\nlso(t-t_2)[t(v)\sech(x)](x-x_0-tv) \\
&+ e^{-itv^2/2}e^{it_2/2}e^{-ixv}\nlso(t-t_2)[r(v)\sech(x)](x+x_0+tv) \\
&+\mathcal{O}(v^{1-\frac{3}{2}\delta}), \qquad t_2\leq t \leq t_3
\end{aligned}
\end{equation}
This is proved by a perturbative argument that enables us to evolve forward a time $(1-\delta)\log v$ at the expense of enlarging the error by a multiplicative factor of $e^{(1-\delta)\log v} = v^{1-\delta}$.  The error thus goes from $v^{-\delta/2}$ at $t=t_2$ to $v^{1-\frac{3}{2}\delta}$ at $t=t_3$.

Now we turn to the details.

\subsection{Phase 1}

Let $u_1(x,t) = \nlso(t)u_0(x)$ and $u(x,t)=\nlsq(t)u_0(x)$.  Let $w= u-u_1$.  Recall that $t_1= |x_0|/v-v^{-\delta}$ so that $x_0+vt_1=-v^{1-\delta}$.  Note that
$$u_1(x,t) = e^{-itv^2/2}e^{it/2}e^{ixv}\sech(x-x_0-tv)$$
We will need the following perturbation lemma.

\begin{lem}
\label{L:approx1}
If $t_a<t_b$, $t_b-t_a \leq c_1$, and $\|w(\cdot,t_a)\|_{L^2}+q\|u_1(0,t)\|_{L_{[t_a,t_b]}^\infty} \leq 1$, then 
$$\|w\|_{L_{[t_a,t_b]}^\infty L_x^2} \leq c_2(\|w(\cdot, t_a)\|_{L_x^2}+q\|u_1(0,t)\|_{L_{[t_a,t_b]}^\infty})$$ 
where the constants $c_1$ and $c_2$ depend only on constants appearing in the Strichartz estimates and are, in particular, independent of $q$ and $v$.
\end{lem}

\begin{proof}
$w$ solves
\begin{align*}
i\partial_t w + \partial_x^2w -q\delta_0(x)w &= -|w+u_1|^2(w+u_1) +|u_1|^2u_1 + q\delta_0(x)u_1\\
&=-\underbrace{|w|^2w}_{\textnormal{cubic}} - \underbrace{(2u_1|w|^2 + \bar u_1 w^2)}_{\textnormal{quadratic}} - \underbrace{(2|u_1|^2w + u_1^2 \bar w)}_{\textnormal{linear}} +q \delta_0(x)u_1\\
\end{align*}
From this equation, $w$ is estimated using Proposition \ref{p:Str}.  For the cubic nonlinear term we take $\tilde p  = \tilde r = 6/5$ and estimate by H\"older as
$$\| |w|^2 w\|_{L_{[t_a,t_b]}^{6/5}L_x^{6/5}} \leq (t_b-t_a)^{1/2}\|w\|_{L_{[t_a,t_b]}^6L_x^6}^2 \|w\|_{L_{[t_a,t_b]}^\infty L_x^2}$$
Since complex conjugates becomes irrelevant in the estimates, both quadratic terms are treated identically.  In Proposition \ref{p:Str}, we take $\tilde p  = \tilde r = 6/5$ and estimate by H\"older as
\begin{align*}
\| u_1 w^2 \|_{L_{[t_a,t_b]}^{6/5}L_x^{6/5}} &\leq (t_b-t_a)^{1/2}\|w\|_{L_{[t_a,t_b]}^6L_x^6}^2 \|u_1\|_{L_{[t_a,t_b]}^\infty L_x^2} \\
&\leq \sqrt 2(t_b-t_a)^{1/2} \|w\|_{L_{[t_a,t_b]}^6L_x^6}^2
\end{align*}
For the linear terms (both of the form $u_1^2w$), we take $\tilde p  = \tilde r = 6/5$ in Proposition \ref{p:Str} and estimate as
\begin{align*}
\| u_1^2 w \|_{L_{[t_a,t_b]}^{6/5}L_x^{6/5}} &\leq (t_b-t_a)^{1/2}\|w\|_{L_{[t_a,t_b]}^6L_x^6} \|u_1\|_{L_{[t_a,t_b]}^6 L_x^6} \|u_1\|_{L_{[t_a,t_b]}^\infty L_x^2} \\
&\leq 2(t_b-t_a)^{2/3} \|w\|_{L_{[t_a,t_b]}^6L_x^6}
\end{align*}
The delta term is estimated by the concluding sentence of Proposition \ref{p:Str} as 
$$q\|u(0,t)\|_{L_{[t_a,t_b]}^{4/3}} \leq q (t_b-t_a)^{3/4} \|u(0,t)\|_{L_{[t_a,t_b]}^\infty}$$
Since $t_b-t_a \leq 1$, collecting the above estimates we have (taking $\|w\|_X = \|w\|_{L_{[t_a,t_b]}^\infty L_x^2} + \|w\|_{L_{[t_a,t_b]}^6L_x^6}$),
$$\|w\|_X \leq c\|w(\cdot, t_a)\|_{L_x^2} + c(t_b-t_a)^{1/2}(\|w\|_X +\|w\|_X^2 + \|w\|_X^3) + cq\|u(0,t)\|_{L_{[t_a,t_b]}^\infty}$$
Provided $(t_b-t_a)^{1/2} \leq 1/(2c)$ above, the linear term on the right can be absorbed by the left as
$$\|w\|_X \leq 2c\|w(\cdot, t_a)\|_{L_x^2} + 2c(t_b-t_a)^{1/2}(\|w\|_X^2 + \|w\|_X^3) + 2cq\|u(0,t)\|_{L_{[t_a,t_b]}^\infty}$$
Continuity of $\|w\|_{X(t_b)}$ as a function of $t_b$ shows that provided $2c(t_b-t_a)^{1/2}(4c\|w(\cdot,t_a)\|_{L^2}+4cq\|u_1(0,t)\|_{L_{[t_a,t_b]}^\infty}) \leq 1/2$, the above estimate implies 
$$\|w\|_X \leq 4c\|w(\cdot, t_a)\|_{L_x^2} + 4cq\|u(0,t)\|_{L_{[t_a,t_b]}^\infty}$$
concluding the proof.
\end{proof}

Now we proceed to apply Lemma \ref{L:approx1}.  The constants $c_1$ and $c_2$ will, for convenience of exposition, be taken to be $c_1=1$ and $c_2=2$.
Let $k\geq 0$ be the integer such that $k \leq t_1 < k+1$.  (Note that $k=0$ if the soliton starts within a distance $v$ of the origin, i.e. $-v-v^{1-\delta}\leq x_0\leq -v^{1-\delta}$, and the inductive analysis below is skipped.) Apply Lemma \ref{L:approx1} with $t_a=0$, $t_b=1$ to obtain (since $w(\cdot,0)=0$)
$$\|w\|_{L_{[0,1]}^\infty L_x^2} \leq 2q\|u_1(0,t)\|_{L_{[0,1]}^\infty}\leq 2q\sech(x_0+v)$$
Apply Lemma \ref{L:approx1} again with $t_a=1$, $t_b=2$ to obtain
\begin{align*}
\|w\|_{L_{[1,2]}^\infty L_x^2} &\leq 2(\|w(\cdot,1)\|_{L_x^2}+ q\|u_1(0,t)\|_{L_{[1,2]}^\infty}) \\
&\leq 2^2q\sech(x_0+v)+2^1q\sech(x_0+2v)
\end{align*}
We continue inductively up to step $k$, and then collect all $k$ estimates to obtain the following bound on the time interval $[0,k]$
$$\|w\|_{L_{[0,k]}^\infty L_x^2} \leq 2q \sum_{j=1}^k 2^{k-j}\sech(x_0+jv)$$
The estimate $\sech \alpha \leq 2e^{-|\alpha|}$ reduces matters to bounding
$$2^kq e^{x_0+v}\sum_{j=0}^{k-1}2^{-j}e^{jv}$$
and, after summing the geometric series, we obtain
$$\|w\|_{L_{[0,k]}^\infty L_x^2} \leq c2^k e^{x_0+v} \frac{(2^{-1}e^{v})^k-1}{2^{-1}e^v-1} \leq cqe^{x_0+kv}$$
where the last inequality requires $2^{-1}e^v \geq 2$.  Finally, applying Lemma \ref{L:approx1} on $[k,t_1]$,
$$\|w\|_{L_{[0,t_1]}^\infty L_x^2} \leq c(qe^{x_0+kv}+q\sech(x_0+t_1v)) \leq cq e^{-v^{1-\delta}}$$
As a consequence, \eqref{E:approx1} follows.

\subsection{Phase 2}

We shall need a lemma stating that the free nonlinear flow is approximated by the free linear flow, and that the perturbed nonlinear flow is approximated by the perturbed linear flow.  Both estimates are consequences of the corresponding Strichartz estimates (Proposition \ref{p:Str}).  Crucially, the hypotheses and estimates of this lemma depend only on the $L^2$ norm of the initial data $\varphi$.  Below, \eqref{E:approxq} is applied with $\varphi(x)=u(x,t_1)$, and $\|u(x,t_1)\|_{L_x^2}=\|u_0\|_{L^2}$ is independent of $v$; thus $v$ does not enter adversely into the analysis.

\begin{lem}
\label{L:approx}
If $\varphi\in L^2$ and $0<t_b$ such that $t_b<c_1\|\varphi\|_{L^2}^{-4}$, then
\begin{equation}
\label{E:approx0}
\|\nlso(t)\varphi - e^{-itH_0}\varphi\|_{L_{[0,t_b]}^\infty L_x^2} \leq c_2t_b^{1/2}\| \varphi\|_{L^2}^3
\end{equation}
\begin{equation}
\label{E:approxq}
\|\nlsq(t)\varphi - e^{-itH_q}\varphi\|_{L_{[0,t_b]}^\infty L_x^2} \leq c_2t_b^{1/2} \| \varphi\|_{L^2}^3
\end{equation}
where $c_1$ and $c_2$ depend only on constants appearing in the Strichartz estimates.  In particular, they are independent of $q$.
\end{lem}

\begin{proof}
Estimate \eqref{E:approx0} is in fact a special case of \eqref{E:approxq} obtained by taking $q=0$.  Let $h(t)=\nlsq(t)\varphi$ so that 
$$i\partial_t h +\tfrac{1}{2}\partial_x^2 h - q\delta_0(x) h + |h|^2h =0$$
with $h(x,0)=\varphi(x)$. Let us define 
\[ X = L_{[0,t_b]}^\infty L_x^2 \cap L_{[0,t_b]}^6L_x^6 \,,\]
with the natural norm, $ \| \bullet \|_X $. 
We apply Proposition \ref{p:Str} with, in the notation of that proposition,  $u(t)=h(t) - e^{-itH_q}\varphi$, $f=-|h|^2h$, $p=r=6$, $\tilde p = \tilde r = 6/5$, and then again with $p=\infty$, $r=2$, $\tilde p = \tilde r = 6/5$, to obtain
\[  \|h(t)  - e^{-itH_q}\varphi\|_X \leq c\| |h|^2 h\|_{L_{[0,t_b]}^{6/5}L_x^{6/5}} 
\,. \]
The generalized H\"older inequality, 
$$ \| h_1 h_2 h_3 \|_p \leq
\| h_1 \|_{q_1} \|h_2\|_{q_2} \|h_3\|_{q_3} \,, \ \ \frac{1}{p} = \frac{1}{
q_1} + \frac{1}{ q_2} + \frac{1 } {q_3} \,, $$
applied with $ h_j = h $, 
$ p = 6/5 $ and $ q_1 = q_2 = 6 $, $ q_3 = 2 $, gives
\begin{equation}
\label{E:110}
\begin{split}
 \|h(t)  - e^{-itH_q}\varphi\|_X 
 &\leq C \|h\|_{L_{[0,t_b]}^6L_x^6}^2 \|h\|_{L^2_{[0,t_b]} L_x^2 }  \\
 &\leq C t_b^{1/2}\|h\|_{L_{[0,t_b]}^6L_x^6}^2 \|h\|_{L_{[0,t_b]}^\infty L_x^2} \\
 &\leq C t_b^{1/2}\|h\|_X^3\,.
\end{split}
\end{equation}
Another application of the homogeneous Strichartz estimate shows that
\[ \| e^{ - i t H_q } \varphi \|_X \leq C \| \varphi \|_{L^2} \,, \]
and consequently,
$$\|h \|_X \leq c  \|\varphi\|_{L^2} +  ct_b^{1/2}\|h\|_X^3$$
By continuity of $\|h\|_{X(t_b)}$ in $t_b$, if $ct_b^{1/2}(2c\|\varphi\|_{L^2})^2 \leq 1/2$,
$$\| h\|_X \leq 2c\|\varphi\|_{L^2}$$
Substituting into \eqref{E:110} yields the result.
\end{proof}

Now we proceed to apply Lemma \ref{L:approx}.  Set $t_2=t_1+2v^{-\delta}$, and apply \eqref{E:approxq} on $[t_1,t_2]$ to obtain
\begin{align*}
u(\cdot,t) &= \nlsq(t-t_1)[u(\cdot,t_1)] \\
&= e^{-i(t-t_1)H_q}[u(\cdot,t_1)] + \mathcal O(v^{-\delta/2})
\end{align*}
By combining this with \eqref{E:approx1},
\begin{equation}
\label{E:approx2}
u(\cdot,t) = e^{-it_1v^2/2}e^{it_1/2}e^{-i(t-t_1)H_q}[e^{ixv}\sech(x-x_0-t_1v)] + \mathcal O(v^{-\delta/2})
\end{equation}

By Proposition \ref{p:as} with $\theta(x)=1$ for $x\leq -1$ and $\theta(x)=0$ for $x\geq 0$, $\varphi(x) = \sech(x)$, and $x_0$ replaced by $x_0+t_1v$,
\begin{equation}
\label{E:approx3}
\begin{aligned}
\indentalign e^{-i(t_2-t_1)H_q}[e^{ixv}\sech(x-x_0-vt_1)](x) \\
&=
\begin{aligned}[t]
&t ( v ) e^{-i(t_2-t_1)H_0}[e^{ixv}\sech(x-x_0-vt_1)](x) \\
&+ r ( v )  e^{-i(t_2-t_1)H_0}[e^{-ixv}\sech(x+x_0+vt_1)](x)  \\
&+ {\mathcal O} (v^{-1}) 
\end{aligned}
\end{aligned}
\end{equation}
By combining \eqref{E:approx2}, \eqref{E:approx3} and \eqref{E:approx0},
$$u(\cdot,t) =
\begin{aligned}[t]
&t(v)e^{-it_1v^2/2}e^{it_1/2}\nlso(t_2-t_1)[e^{ixv}\sech(x-x_0-vt_1)](x) \\
&+r(v)e^{-it_1v^2/2}e^{it_1/2}\nlso(t_2-t_1)[e^{-ixv}\sech(x+x_0+vt_1)](x)\\
&+\mathcal{O}(v^{-\delta/2})
\end{aligned}
$$
By noting that
\begin{align*}
\indentalign \nlso(t_2-t_1)[e^{ixv}\sech(x-x_0-t_1v)] \\
&= e^{-i(t_2-t_1)v^2/2}e^{i(t_2-t_1)/2}e^{ixv}\sech(x-x_0-t_2v)
\end{align*}
and
\begin{align*}
\indentalign \nlso(t_2-t_1)[e^{-ixv}\sech(x+x_0+t_1v)] \\
&= e^{-i(t_2-t_1)v^2/2}e^{i(t_2-t_1)/2}e^{-ixv}\sech(x+x_0+t_2v)
\end{align*}
we obtain \eqref{E:approx4}.

\subsection{Phase 3}

Let $t_3=t_2+(1-\delta) \log v$.  Label
$$u_\trans(x,t) =e^{-itv^2/2}e^{it_2/2}e^{ixv}\nlso(t-t_2)[t(v)\sech(x)](x-x_0-tv)$$
for the transmitted (right-traveling) component and
$$u_\refl(x,t) = e^{-itv^2/2}e^{it_2/2}e^{-ixv}\nlso(t-t_2)[r(v)\sech(x)](x+x_0+tv)$$
for the reflected (left-traveling) component.  By Appendix \ref{localization}, for each $k\in \mathbb{N}$, there is a constant $c(k)>0$ and an exponent $\sigma(k)>0$ such that
\begin{equation}
\label{E:trans}
\|u_\trans(x,t)\|_{L_{x<0}^2}+ \|u_\refl(x,t)\|_{L_{x>0}^2} + |u_{\tr}(0,t)| + |u_{\refl}(0,t)| \leq \frac{c(k)(\log v)^{\sigma(k)}}{v^{k(1-\delta)}}
\end{equation}
uniformly on the time interval $[t_2,t_3]$.  We shall need the following perturbation lemma, again a consequence of the Strichartz estimates.

\begin{lem}
\label{L:approx6}
  Let $w=u-u_\trans-u_\refl$.
If $t_a<t_b$, $t_b-t_a \leq c_1$, and 
$$\|w(\cdot,t_a)\|_{L_x^2} + \frac{c(k)\la q\ra(\log v)^{\sigma(k)}}{v^{k(1-\delta)}} \leq 1
\,, $$ 
 then
\begin{align*}
\|w\|_{L_{[t_a,t_b]}^\infty L_x^2} &\leq c_2\left(\|w(\cdot,t_a)\|_{L_x^2} + \frac{c(k)\la q\ra(\log v)^{\sigma(k)}}{v^{k(1-\delta)}}\right)
\end{align*}
The constants $c_1$, $c_2$ depend only on constants appearing in the Strichartz estimates and are in particular independent of $q$ and $v$.
\end{lem}

\begin{proof}
We write the equation satisfied by $w$:
\begin{align*}
\indentalign i\partial_t w + \tfrac{1}{2}\partial_x^2 w - q\delta_0(x)w \\
&= 
\begin{aligned}[t]
&-|w+u_\trans+u_\refl|^2(w+u_\trans+u_\refl) + |u_\trans|^2u_\trans + |u_\refl|^2u_\refl\\
&+q\delta_0(x)u_\trans -q\delta_0(x)u_\refl \\
\end{aligned}\\
&=
\begin{aligned}[t]
&-|w|^2w - (2(u_\trans+u_\refl)|w|^2 + (\bar u_\trans + \bar u_\refl)w^2) - (2|u_\trans+u_\refl|^2w + (u_\trans+u_\refl)^2\bar w)\\
&-\underbrace{(u_\trans^2\bar u_\refl + 2u_\refl|u_\trans|^2 +u_\refl^2\bar u_\trans + 2u_\trans|u_\refl|^2)}_{\trans-\refl\textnormal{ interaction}}+\underbrace{q\delta_0(x)u_\trans}_{\trans-\textnormal{delta}} -\underbrace{q\delta_0(x)u_\refl}_{\refl-\textnormal{delta}} \\
\end{aligned}
\end{align*}
$w$ is estimated using Proposition \ref{p:Str} and Proposition \ref{p:Str}.  The cubic, quadratic, and linear in $w$ terms on the first line are estimated exactly as was done in the proof of Lemma \ref{L:approx1}.  For the ``$\trans-\refl$ interaction terms'' (taking $u_\refl|u_\trans|^2$ as a representative example), we apply Proposition \ref{p:Str} with $\tilde p = 4/3$, $\tilde r = 1$ and estimate as
\begin{equation}
\label{E:120}
\| u_\refl |u_\trans|^2 \|_{L_{[t_a,t_b]}^{4/3}L_x^1} \leq c(t_b-t_a)^{3/4}\|u_\trans\|_{L_{[t_a,t_b]}^\infty L_x^2}\|u_\refl u_\trans\|_{L_{[t_a,t_b]}^\infty L_x^2}
\end{equation}
$\|u_\trans\|_{L_x^2} = \sqrt 2 |t(v)|$ by mass conservation for the free nonlinear flow, and
\begin{align*}
\|u_\trans u_\refl\|_{L_{[t_a,t_b]}^\infty L_x^2} &\leq \|u_\trans u_\refl\|_{L_{[t_a,t_b]}^\infty L_{x<0}^2} + \|u_\trans u_\refl\|_{L_{[t_a,t_b]}^\infty L_{x>0}^2}\\
&\leq \|u_\refl\|_{L_{[t_a,t_b]}^\infty L_x^\infty} \|u_\trans\|_{L_{[t_a,t_b]}^\infty L_{x<0}^2} + \|u_\trans\|_{L_{[t_a,t_b]}^\infty L_x^\infty} \|u_\refl\|_{L_{[t_a,t_b]}^\infty L_{x>0}^2}
\end{align*}
Now
\begin{align*}
\|u_\refl\|_{L_t^\infty L_x^\infty} &= \|\nlso(t)[r(v)\sech](x)\|_{L_t^\infty L_x^\infty} \\
&\leq \|\nlso(t)[r(v)\sech](x)\|_{L_t^\infty L_x^2}^{1/2} \|\partial_x \nlso(t)[r(v)\sech](x)\|_{L_t^\infty L_x^2}^{1/2} \\
&\leq c
\end{align*}
by mass and energy conservation of the free nonlinear flow.  Similarly, $\|u_\trans\|_{L_t^\infty L_x^\infty} \leq c$.  By this and \eqref{E:trans}, the above yields
$$\|u_\trans u_\refl\|_{L_{[t_a,t_b]}^\infty L_x^2} \leq \frac{ c(k)(\log v)^{\sigma(k)}}{v^{k(1-\delta)}}$$
Thus, by \eqref{E:120},
\begin{equation}
\label{E:121}
\| u_\refl |u_\trans|^2 \|_{L_{[t_a,t_b]}^{4/3}L_x^1} \leq \frac{c(k)(\log v)^{\sigma(k)}}{v^{k(1-\delta)}}
\end{equation}
and similarly for all other ``$\trans-\refl$ interaction'' terms.  Now we address the ``$\trans-$delta'' and ``$\refl-$delta'' terms (working with $q\delta_0(x)u_\trans$ as the representative of both).  By Proposition \ref{p:Str}, we estimate as
$$q \|u_\trans(0,t)\|_{L_{[t_a,t_b]}^{4/3}} \leq c(t_b-t_a)^{3/4}q \|u_\trans(0,t)\|_{L_{[t_a,t_b]}^\infty}$$
By \eqref{E:trans},
\begin{equation}
\label{E:122}
q \|u_\trans(0,t)\|_{L_{[t_a,t_b]}^{4/3}} \leq c(k)q(t_b-t_a)^{3/4} \frac{(\log v)^{\sigma(k)}}{v^{k(1-\delta)}}
\end{equation}
Collecting \eqref{E:121}, \eqref{E:122}, and the estimates for cubic, quadratic, and linear terms in $w$ (as exposed in Lemma \ref{L:approx1}), we have, with $\|w\|_X = \|w\|_{L_{[t_a,t_b]}^\infty L_x^2} + \|w\|_{L_{[t_a,t_b]}^6 L_x^6}$
$$\|w\|_X \leq c\|w(\cdot, t_a)\|_{L^2} + c(t_b-t_a)^{1/2}(\|w\|_X+\|w\|_X^2+\|w\|_X^3) + \frac{c(k)\la q\ra(\log v)^{\sigma(k)}}{v^{k(1-\delta)}}$$
If $c(t_b-t_a)^{1/2}\leq \frac{1}{2}$, then the first-order $w$-term on the right side can be absorbed by the left, giving
$$\|w\|_X \leq 2c\|w(\cdot, t_a)\|_{L^2} + 2c(t_b-t_a)^{1/2}(\|w\|_X^2+\|w\|_X^3) + \frac{2c(k)\la q\ra(\log v)^{\sigma(k)}}{v^{k(1-\delta)}}$$
By continuity of $\|w\|_{X(t_b)}$ in $t_b$, if 
$$2c(t_b-t_a)^{1/2}\left( 4c\|w(\cdot,t_a)\|_{L^2} + \frac{4c(k)\la q\ra (\log v)^{\sigma(k)}}{v^{k(1-\delta)}} \right) \leq \frac{1}{2}$$
we have
$$\|w\|_X \leq 4c\|w(\cdot, t_a)\|_{L^2} + \frac{4c(k) \la q\ra(\log v)^{\sigma(k)}}{v^{k(1-\delta)}}$$
completing the proof.
\end{proof}

Assume that $\alpha = q/v$ has been fixed.  Choose $k=k(\delta)$ large so that $k(1-\delta)\geq 3$.  Then the coefficient appearing in Lemma \ref{L:approx6} is bounded by
$$\frac{c(k)\la q \ra(\log v)^{\sigma(k)}}{v^{k(1-\delta)}} \leq c(k) \left< \frac{q}{v}\right> \frac{(\log v)^{\sigma(k)}}{v^2}$$
Now take $v$ sufficiently large in terms of $\la q/v \ra$ and $k$ (thus in terms of $\delta$) so that the above is bounded by $v^{-1}$.

Now we implement Lemma \ref{L:approx6}.  For convenience of exposition, we take $c_1=1$, $c_2=2$.  Let $\ell$ be the integer such that $\ell <(1-\delta)\log v < \ell+1$.  We then apply Lemma \ref{L:approx6} successively on the intervals $[t_2,t_2+1], \ldots, [t_2+\ell-1,t_2+\ell]$ as follows.  Applying Lemma \ref{L:approx6} on $[t_2,t_2+1]$, we obtain
$$\|w(\cdot,t)\|_{L_{[t_2,t_2+1]}^\infty L_x^2} \leq 2(\|w(\cdot,t_2)\|_{L_x^2} + v^{-1})$$
Applying Lemma \ref{L:approx6} on $[t_2+1,t_2+2]$ and combining with the above estimate,
$$\|w(\cdot,t)\|_{L_{[t_2+1,t_2+2]}^\infty L_x^2} \leq 2^2\|w(\cdot,t_2)\|_{L_x^2} + (2^2+2)v^{-1}$$
Continuing up to the $\ell$-th step and then collecting all of the above estimates,$$\|w(\cdot,t)\|_{L_{[t_2,t_3]}^\infty L_x^2} \leq 2^\ell\|w(\cdot,t_2)\|_{L_x^2} + (2^\ell + \cdots +2)v^{-1} $$
Since $\|w(\cdot, t_2)\|_{L_x^2} \leq v^{-\delta/2}$ and $2^\ell \leq v^{1-\delta}$, 
\begin{equation}
\label{E:approx7}
\|w(\cdot,t)\|_{L_{[t_2,t_3]}^\infty L_x^2} \leq cv^{1-\frac{3}{2}\delta}
\end{equation}
thus proving \eqref{E:post}. 
 
Now we complete the proof of the main theorem and obtain \eqref{eq:th}.  By \eqref{E:approx7} and \eqref{E:trans},
\begin{equation}
\label{E:approx8}
\|u(\cdot,t)-u_\trans(\cdot,t)\|_{L_{x>0}^2} \leq \|w(\cdot,t)\|_{L_{x>0}^2}  + \|u_\refl(\cdot,t)\|_{L_{x>0}^2} \leq cv^{1-\frac{3}{2}\delta}
\end{equation}
Since $\|u_\trans(\cdot,t)\|_{L_x^2} = t(v)$, \eqref{E:trans} implies $\|u_\trans(\cdot,t)\|_{L_{x>0}^2} = t(v) + \mathcal{O}(v^{1-\frac{3}{2}\delta})$, which combined with \eqref{E:approx8} gives \eqref{eq:th} and proves Theorem \ref{th:1}.

\section{Resolution of outgoing waves}
\label{resolution}

In this section, we prove Theorem \ref{th:2}.  We note that the proof of Theorem \ref{th:1} presented in \S\ref{proof} in fact provided a more complete long-time description of the solution:
\begin{equation}
\label{E:asymp}
u(x,t) = 
\begin{aligned}[t]
& e^{-itv^2/2}e^{it_2/2}e^{ixv}\nlso(t-t_2)[t(v)\sech](x-x_0-tv) \\
&+ e^{-itv^2/2}e^{it_2/2}e^{-ixv}\nlso(t-t_2)[r(v)\sech](x+x_0+tv) \\
&+\mathcal{O}_{L_x^2}(v^{1-\frac{3}{2}\delta}), \qquad \frac{|x_0|}{v} +v^{-\delta}\leq t \leq c(1-\delta)\log v
\end{aligned}
\end{equation}
where $t(v)$, $r(v)$ are defined in \eqref{eq:tr} and $\nlso(t)\varphi$ denotes the solution to the NLS equation $i\partial_th + \tfrac{1}{2}\partial_x^2h + |h|^2h =0$ (without potential) and initial data $h(x,0) = \varphi(x)$.  It thus suffices to obtain the resolution of $\nlso(t-t_2)[t(v)\sech]$ and $\nlso(t-t_2)[r(v)\sech]$ into solitons plus radiation decaying in $L_x^\infty$.  By the phase invariance of the free nonlinear flow
$$ \nlso(t-t_2)[t(v)\sech] = \frac{t(v)}{|t(v)|}\nlso(t-t_2)[|t(v)|\sech]$$
and similarly for $\nlso(t-t_2)[r(v)\sech]$.
Since $0\leq |t(v)|, |r(v)| \leq 1$, we apply asymptotics \eqref{eq:apb} proved of Appendix \ref{alphasech} using the  inverse scattering method. When 
$ | t ( v) | $ or $ | r ( v ) | $ is equal to $ 1/2 $ we use the result
of \cite{Ka} recalled in \eqref{eq:apb1}.  The result obtained by these substitutions differs from that stated in Theorem \ref{th:2} by a factor of 
\begin{equation}
\label{E:deltashift}
\exp \left(i \, \frac{1-A_T^2}{2} \cdot v^{-\delta}\right)
\end{equation}
for $u_T(x,t)$, owing to the fact that $t_2=|x_0|/v+ v^{-\delta}$.  But \eqref{E:deltashift} differs from $1$ by $\sim v^{-\delta}$, and thus omitting it only introduces a discrepancy of $v^{-\delta}$ in both $L_x^2$ and $L_x^\infty$.  There is a similar inconsequential disparity in the $u_R(x,t)$ part.

\appendix
\section{Spatial localization of the free nonlinear propagation}
\label{localization}

Let $\varphi \in \mathcal{S}$ and 
\begin{equation}
\left\{
\begin{aligned}
&i \partial_t h + \tfrac{1}{2} \partial_x^2 h + |h|^2 h  =  0, \\
&h(x,0)  =  \varphi(x).
\end{aligned}
\right.
\end{equation}
\textbf{Notational conventions}.  We denote $\partial_x$ by $\partial$ hereafter.  The $x$ and $t$ dependence of $h(x,t)$ will be routinely dropped.  The constants $c(k)$, $\sigma(k)$, and the polynomials $g_\gamma(t)$ that appear below may change (enlarge) from one line to the next without comment. The constants $ c( k ) $ depend on the fixed function $ \varphi \in {\mathcal S} $.

The solution $h$ satisfies conservation of mass and conservation of energy, which means that the integrals
$$ E_0 = \int_\RR  |h|^2 dx \,, \ \ 
 E_2 = - \int_\RR  (| \partial h |^2 - |h|^4) dx \,, $$
are independent of time $ t$.  Since $\|h\|_{L^\infty}^2 \leq \|h\|_{L^2}\|\partial h\|_{L^2}$, we have $\|h\|_{L^4}^4 \leq \|h\|_{L^2}^3\|\partial h\|_{L^2}$ and it follows from the $E_2$ and $E_0$ conservation that $\|\partial h\|_{L^2} \leq c$, where $c$ depends on $\|\varphi\|_{L^2}$ and $\|\partial \varphi \|_{L^2}$.

In fact, there are an infinite number of conserved integrals, $ E_k $, with 
integrands defined inductively as follows
\begin{equation}
\label{E:A3}
f_0  = |h|^2  \,, \ \ \
 f_{k+1} =  h \partial \left( \frac{1}{h} f_k \right) + \sum_{j_1+j_2=k-1} 
f_{j_1} f_{j_2}, 
\end{equation}
see \cite[\S 8]{ZS72} for a proof of this fact (rescaling time and putting 
$ \kappa = 2 $ produces an agreement with our slightly different 
convention).  The inductive definition of $f_k$ and the Sobolev 
embedding theorem can now be used to show that, for $\ell\geq 2$, 
\begin{equation}
\label{E:A4}
 E_{2 \ell}  = ( -1)^{\ell} \int_{\RR } | \partial^\ell h ( x) |^2 d x 
+ {\mathcal O} ( ( 1 + \| h \|_{H^{\ell - 1}_x } ) ^{2\ell + 2} ) \,, 
\end{equation}
and hence for $\ell\geq 0$, we have
\begin{equation}
\label{E:ZS}
\|\partial^\ell h \|_{L^2} \leq c(\ell)
\end{equation}
where $c(\ell)$ depends upon Sobolev norms of the initial data $\varphi$ of at most order $\ell$.
We now elaborate on how to obtain \eqref{E:A4}.  An inductive argument using \eqref{E:A3} shows that for $k\geq 0$, $f_k$ is of the form
\begin{equation}
\label{E:A6}
f_k = h\partial^k \bar h + h \sum_{j\geq 1,\; 2j\leq k} p(2j+1,k-2j)
\end{equation}
where $p(n,m)$ indicates a linear combination of terms of degree $n$ and \textit{cumulative} order $m$, or more precisely terms of the form
\begin{equation}
\label{E:A5}
\partial^{\alpha_1}\tilde h  \, \partial^{\alpha_2}\tilde h\cdots \partial^{\alpha_n}\tilde h, \quad \alpha_1+\cdots+\alpha_n=m
\end{equation}
and $\tilde h$ is either $h$ or $\bar h$.  To prove \eqref{E:A4} for $\ell\geq 2$, one uses \eqref{E:A6} for $k=2\ell$ and it only remains to verify that for any $n\geq 4$ and $m\leq 2\ell-2$,
\begin{equation}
\label{E:A7}
\int  p(n,m) \, dx \leq \|h\|_{H^{\ell-1}}^n
\end{equation}
We now show this.   Note that in \eqref{E:A5}, we may assume without loss of generality that $\alpha_1\leq \alpha_2 \leq \cdots \leq \alpha_n$.  

\noindent\textbf{Case 1}.  $\alpha_n\leq \ell-1$.  It follows that $\alpha_j \leq \ell-2$ for all $j\leq n-2$ and $\alpha_{n-1}\leq \ell-1$.  We estimate as:
$$\left| \int \partial^{\alpha_1}\tilde h  \, \partial^{\alpha_2}\tilde h\cdots \partial^{\alpha_n}\tilde h \, dx\right| \leq \left( \prod_{j=1}^{n-2} \|\partial^{\alpha_j} h\|_{L^\infty} \right) \|\partial^{\alpha_{n-1}}h\|_{L^2} \|\partial^{\alpha_n}h\|_{L^2}\leq c\|h\|_{H^{\ell-1}}^n$$
by Sobolev embedding.

\noindent\textbf{Case 2}.  $\alpha_n\geq \ell$.  In this case, we begin by integrating by parts to obtain
\begin{equation}
\label{E:A8}
(-1)^{\alpha_n-\ell+1}\int \partial^{\alpha_n-\ell+1}( \partial^{\alpha_1}\tilde h  \cdots \partial^{\alpha_{n-1}}\tilde h) \; \partial^{\ell-1}\tilde h \, dx
\end{equation}
The Leibniz rule expansion is
\begin{equation}
\label{E:A9}
\partial^{\alpha_n-\ell+1}( \partial^{\alpha_1}\tilde h  \cdots \partial^{\alpha_{n-1}}\tilde h) = \sum c_\mu \partial^{\mu_1+\alpha_1}\tilde h \cdots \partial^{\mu_{n-1}+\alpha_{n-1}}\tilde h
\end{equation}
where the sum is over $(n-1)$-tuples $\mu$ such that $\mu_1+\cdots+\mu_{n-1} = \alpha_n-\ell+1$ and $c_\mu$ is some constant depending on $\mu$.  By adding the $\alpha$ and $\mu$ constraints, we obtain that $(\mu_1+\alpha_1)+\cdots +(\mu_{n-1}+\alpha_{n-1}) \leq \ell-1$ and thus there is at most one index $j_*$  ($1\leq j_* \leq n-1$) such that $\mu_{j_*}+\alpha_{j_*} = \ell-1$ and for all remaining $j$ ($1\leq j\leq n-1$, $j\neq j_*$) we have $\mu_j+\alpha_j \leq \ell-2$.  (If no such $j_*$ exists, take $j_*$ to be any fixed index $1\leq j_*\leq n-1$.)   By substituting \eqref{E:A9} into \eqref{E:A8}, we estimate as
$$\left|\int \partial^{\alpha_1}\tilde h  \, \partial^{\alpha_2}\tilde h\cdots \partial^{\alpha_n}\tilde h \, dx \right| \leq \left( \prod_{j=1, j\neq j_*}^{n-1} \|\partial^{\mu_j+\alpha_j} h\|_{L^\infty} \right) \|\partial^{\mu_{j_*}+\alpha_{j_*}}h\|_{L^2}\|\partial^{\ell-1}h\|_{L^2} \leq c\|h\|_{H^{\ell-1}}^n$$
again by Sobolev embedding.  This concludes the proof of \eqref{E:A7}, thus \eqref{E:A4}, and thus \eqref{E:ZS}.

Using that the commutator $[(x+it\partial),i\partial_t+\tfrac{1}{2}\partial^2]=0$  and some integration by parts manipulations, we have the pseudoconformal conservation law:
$$
\int_x |(x+it\partial) h(x,t)|^2 dx - t^2 \int_x |h(x,t)|^4 dx + \int_0^t s \int_x |h(x,s) |^4 dx ds = \int_x |x \varphi(x)|^2 dx.
$$
From this, \eqref{E:ZS} for $\ell=0,1$, and the Gagliardo-Nirenberg estimate $\|h\|_{L^4}^4 \leq \|h\|_{L^2}^3\|\partial h\|_{L^2}$, we have
$$ \| x h \|_{L^2}  \leq  c \langle t \rangle$$
where $c$ depends on $\|x\varphi\|_{L^2}$, $\|\varphi\|_{L^2}$, and $\|\partial \varphi\|_{L^2}$.  We want to show that more generally, for each $k\in \mathbb{Z}$, $k\geq 0$, we have
\begin{equation}
\label{E:kthbound}
\| x^\alpha \partial^\beta h \|_{L^2} \leq c(k) \la t \ra^{\sigma(k)} \qquad \text{for }\alpha+\beta =k, \; \alpha,\beta \geq 0, \; \alpha,\beta \in \mathbb{Z}
\end{equation}
Here $c(k)$ is a constant depending on $k$ and weighted Sobolev norms of the initial data $\varphi$ (up to order $2k$), and $\sigma(k)$ is a positive exponent depending upon $k$.  We are not concerned with obtaining the optimal value of $\sigma(k)$; the mere fact that the bound in \eqref{E:kthbound} is power-like in $t$, as opposed to exponential in $t$, suffices for our purposes.  In our proof, both $c(k)$ and $\sigma(k)$ will be increasing with $k$, and will go to $+\infty$ as $k\to +\infty$.

Let $\Lambda_0 = \partial$ and $\Lambda_1 = (x+it\partial)$.  Note that both operators have the commutator property
\begin{equation}
\label{E:comm}
[\Lambda_j, (i\partial_t + \tfrac12\partial^2)]=0\,, \ \ j=0, 1\,. 
\end{equation}
We first claim that for each $k\geq 0$,  there exists a constant $c(k)>0$ and an exponent $\sigma(k)>0$ such that 
\begin{equation}
\label{E:A1}
\| \Lambda_{j_1} \Lambda_{j_2} \cdots \Lambda_{j_k} h \|_{L^2} \leq c(k)\la t \ra^{\sigma(k)} \text{ for all }j_1, \ldots j_k \in \{ 0, 1\} \,. 
\end{equation}
 When we wish to consider a composition of the form $\Lambda_{j_1} \Lambda_{j_2} \cdots \Lambda_{j_k}$ and do not care to report whether each operator in the composition is $\Lambda_0$ or $\Lambda_1$, we will instead write the composition as $\Lambda^k$.  We prove \eqref{E:A1} by induction on $k$ .  When $k=0$, \eqref{E:A1} is just the mass conservation law.  Suppose that \eqref{E:A1} holds for $0, \ldots, k-1$; we aim to prove it holds for $k$.  The main ingredient (in addition to the inductive hypothesis) is \eqref{E:ZS}.  Fix $j_1, \ldots j_k \in \{0,1\}$, and apply the operator $\Lambda_{j_1}\cdots \Lambda_{j_k}$ to the equation, pair with $-i \overline{ \Lambda_{j_1}\cdots \Lambda_{j_k} h}$, integrate in $x$, take twice the real part, and appeal to \eqref{E:comm} to obtain
\begin{equation}
\label{E:A2}
\partial_t \| \Lambda_{j_1} \cdots \Lambda_{j_k} h \|_{L^2}^2 = 2\Re i\int  \Lambda_{j_1} \cdots \Lambda_{j_k} |h|^2h \;\;  \overline{\Lambda_{j_1}\cdots \Lambda_{j_k} h} \, dx
\end{equation}
Note that
$$\Lambda_0 \, F(h,\bar h)h = \partial F(h,\bar h) \; h + F(h,\bar h) \; \Lambda_0 h$$
and
$$\Lambda_1 \, F(h,\bar h)h = it\partial F(h,\bar h) \; h + F(h,\bar h) \; \Lambda_1 h$$
Both of these product rules take the form
$$\Lambda \, |h|^2h = g(t)\partial F(h,\bar h) \; h + F(h,\bar h) \; \Lambda h$$
where $g(t)$ is a polynomial in $t$ of degree $\leq 1$.  Thus we see that
$$\Lambda_{j_1} \cdots \Lambda_{j_k} \, |h|^2h = |h|^2 \Lambda_{j_1} \cdots \Lambda_{j_k}h + \sum_{\substack{\gamma_1+\gamma_2 =k \\ \gamma_2 \leq k-1}} g_\gamma(t) \partial^{\gamma_1}|h|^2 \; \Lambda^{\gamma_2}h$$
where $g_{\gamma}(t)$ is a polynomial in $t$.  Substituting into \eqref{E:A2}, we obtain two terms: the first is zero since it is the real part of a purely imaginary number; the second is estimated by the H\"older inequality to obtain:
$$| \, \partial_t \| \Lambda_{j_1} \cdots \Lambda_{j_k} h \|_{L^2}^2 \, | \leq c(k) \la t \ra^{\sigma(k)} \left( \sup_{j\leq k-1} \|\partial^j |h|^2 \|_{L^\infty}\right) \left(\sup_{j\leq k-1}\| \Lambda^jh\|_{L^2}\right)\|\Lambda_{j_1} \cdots \Lambda_{j_k} h \|_{L^2} $$
By Sobolev embedding estimates, \eqref{E:ZS}, and the induction hypothesis, we have
$$| \, \partial_t \| \Lambda_{j_1} \cdots \Lambda_{j_k} h \|_{L^2}^2 \, | \leq c(k) \la t \ra^{\sigma(k)} \|\Lambda_{j_1} \cdots \Lambda_{j_k} h \|_{L^2}$$
from which \eqref{E:A1} follows.

Now to deduce \eqref{E:kthbound} from \eqref{E:A1}, we just note that since $x=\Lambda_1-it\Lambda_0$, there are polynomials $g_j(t)$ such that the following relation holds:
$$x^\alpha \partial^\beta = \sum_{j\in\{0,1\}^{\alpha+\beta}} g_j(t)\Lambda_{j_1}\cdots \Lambda_{j_{\alpha+\beta}}$$

Let us now consider the application of \eqref{E:kthbound} to obtain \eqref{E:trans} in the Phase 3 analysis.  We have $x_0+tv \geq v^{1-\delta}$ for $t\geq t_2$.  If $x<0$, then
$$v^{k(1-\delta)} \leq (x_0+tv)^k \leq |x-x_0-tv|^k$$
Thus
\begin{align*}
v^{k(1-\delta)}\|u_{\tr}(x,t)\|_{L^2_{x<0}} &\leq  \|(x-x_0-tv)^k \nlso(t-t_2)[t(v)\sech](x-x_0-tv)\|_{L_x^2} \\
&= \|x^k \nlso(t-t_2)[t(v)\sech](x)\|_{L_x^2} \\
&\leq c(k)\la t-t_2 \ra^{\sigma(k)}
\end{align*}
by \eqref{E:kthbound}, which gives the first estimate in \eqref{E:trans}.  The second is obtained similarly.  To obtain the third, we note that for $t\geq t_2$,\begin{align*}
|u_{\tr} (0,t)|^2 &= |\nlso(t-t_2)[t(v)\sech](-x_0-tv)|^2\\
&= -\int_{-\infty}^0 \partial_x |\nlso(t-t_2)[t(v)\sech](x-x_0-tv)|^2 \, dx
\end{align*}
and this can be estimated by 
\[  \| \nlso(t-t_2)[t(v)\sech](x-x_0-tv) \|_{L_{x<0}^2} \| \partial_x \nlso(t-t_2)[t(v)\sech](x-x_0-tv) \|_{L_{x<0}^2}  \,. 
\]
Using \eqref{E:kthbound} as before establishes
$$v^{k(1-\delta)} \| \partial_x^j 
\nlso(t-t_2)[t(v)\sech] (x-x_0-tv) \|_{L_{x<0}^2} 
\leq c(k+j) \la t-t_2 \ra^{\sigma(k+j)}$$
Replacing $(c(k)c(k+1))^{1/2}$ by $c(k)$ and $\frac{1}{2}(\sigma(k)+\sigma(k+1))$ by $\sigma(k)$, we obtain the bound \eqref{E:trans}.  
Finally, we note that the fourth bound in \eqref{E:trans} is similar to the third.

\section{Free nonlinear evolution of $ \alpha \, \sech $}
\label{alphasech}

This appendix is devoted to showing that for $ 0 < \alpha < 1 $
\begin{equation}
\label{eq:apb}
\nlso ( \alpha \, \sech ) = \left\{ \begin{array}{ll} 
\nlso ( ( 2 \alpha - 1 ) \sech ( ( 2 \alpha - 1)
\bullet ) ) + {\mathcal O}_{L^\infty } ( t^{-\frac12} ) 
& 1/2 < \alpha < 1 \,, \\
{\mathcal O}_{L^\infty } ( t^{-\frac12} ) &  0 < \alpha < 1/2 \,. 
\end{array} \right.
\end{equation}

A more precise understanding of error terms is possible thanks
to advances in the study asymptotics for integrable nonlinear waves
\cite{DIZ}, \cite{Ka}. Since we do not know an exact reference
for \eqref{eq:apb} we present a proof of this 
simpler asymptotic result. It is based on the now classical work 
on the inverse scattering method initiated for NLS by 
Zakharov-Shabat \cite{ZS72} -- see \cite{DIZ},\cite{FT}
for discussion and references. For the reader's convenience, 
especially in view of different conventions used in different
sources for our argument, we review all the needed aspects of the method.

In the case of $ \alpha = 1/2 $
we can use the result of \cite{Ka} to conclude that 
\begin{equation}
\label{eq:apb1} 
\nlso ( \alpha \, \sech ) =
{\mathcal O}_{L^\infty } ( ( \log t /t) ^{\frac12} )\,. 
\end{equation}
A slightly inaccurate 
statement similar to \eqref{eq:apb} was given in \cite{M81} and
the calculation of the scattering matrix in that paper was our
starting point in obtaining \eqref{eq:apb}.

\subsection{Inverse scattering method}
We present a quick review of this celebrated method.
Thus, let us consider two operators acting on $ {\mathcal S}' ( 
\RR ; \CC^2 ) $:
$$  L = -iJ\partial_x + iJ Q \,, \ \
  A = J\partial_x^2 - \frac 12 JQ_x -JQ\partial_x - \frac 12 Q^2J \,, $$ 
where
$$Q = Q (t,  x ) = \begin{bmatrix} \ \ 0 & u ( t, x )  
\\ - \overline {u  ( t, x )}
  &  \ 0 \end{bmatrix}\,, \ \ 
u ( t, \bullet ) \in {\mathcal S} ( \RR ) \,, 
\qquad J = \begin{bmatrix} -1 & 0 \\ 0 & 1\end{bmatrix}\,. $$
Then 
$$[  L,   A] = - \frac i 2 Q_{xx} + iQ^3$$
which is checked by using
$$
JQJ = -Q\,, \ \ J^2=I \,, \ \ 
Q^2 = \begin{bmatrix} -|u|^2 & 0 \\ 0 & -|u|^2 \end{bmatrix}\,. $$
It is now the case that 
\begin{equation}
\label{eq:evolb} 
 \partial_t   L = i[  L,  A]  \ 
\Longleftrightarrow \ i\partial_t u + \frac12 \partial_x^2 u + |u|^2u=0 \,,
\end{equation} 
and, since we are solving NLS,
we assume that these equivalent equations hold.


We now consider scattering theory for the time dependent operator
$ L $. For that we introduce
special solutions to $ L \psi = \lambda \psi $ with prescribed
asymptotic behaviour:
\begin{equation}
\begin{split}
\label{eq:Jost} 
& \bar \psi ( x , \lambda ) \simeq \left[ \begin{array}{l}
e^{-i x \lambda } \\ \  0 \end{array} \right] \,, \ \ 
\psi ( x , \lambda ) \simeq \left[ \begin{array}{l}
 \ 0 \\    e^{i x \lambda }  \end{array} \right] \,, \ \ 
x \longrightarrow + \infty  \\
& \varphi ( x , \lambda ) \simeq \left[ \begin{array}{l}
 e^{-i x \lambda } \\  \  0 \end{array} \right] \,, \ \ 
\bar \varphi ( x , \lambda ) \simeq \left[ \begin{array}{l}
\ 0 \\  e^{i x \lambda }  \end{array} \right]\,,  \ \ 
x \longrightarrow -  \infty \,,
\end{split}
\end{equation}
see instance \cite[Sect.I.5]{FT}. 
Here for vector valued functions, 
$ \bar \varphi \defeq [ \bar \varphi_2 , - \bar \varphi_1 ]^t $, if
$ \varphi = [ \varphi_1 , \varphi_2 ]^t $.
Each pair
of solutions forms a basis for the solution set and, for $ \lambda 
\in \RR $, 
\begin{equation}
\label{eq:relb}
\begin{split}
& \varphi ( x, \lambda ) = a ( \lambda ) \bar \psi ( x , \lambda ) 
+ b ( \lambda ) \psi  ( x , \lambda) \,, \\
& \bar \varphi  ( x, \lambda ) = \bar a ( \lambda )  \psi ( x , \lambda ) 
- \bar b ( \lambda ) \bar \psi ( x , \lambda) \,, \\
&   \ \ \ \ \ \ \ \ \ \ | a ( \lambda ) |^2 + | b ( \lambda ) |^2 = 1 \,. 
\end{split} \end{equation}
Another consequence comes from \eqref{eq:evolb}. If $ L( t ) \psi ( t ) 
= \lambda \psi ( t) $ then we see that 
\[ ( L - \lambda ) ( i \partial_t \psi - A \psi ) = 0 \]
and hence
\[  i \partial_t \psi ( t )  - A \psi ( t ) =  c_1 ( t ) \psi ( t ) 
+ c_2 \bar \psi ( t ) \,.\]
Now we note that for $ u ( t , \bullet ) \in {\mathcal S} 
( \RR ) $, $ A \simeq J \partial_x^2 $, as $ |x| \rightarrow \infty $, 
and  the asymptotic behaviour  \eqref{eq:Jost} 
gives $ c_1 ( t ) \equiv \lambda^2 $, $ c_2 ( t ) \equiv 0 $. 
More generally we conclude that 
\begin{equation}
\label{eq:moreg} \begin{split}
 & i \partial_t \psi = ( A + \lambda^2 ) \psi \,, \ \
i \partial_t \bar \psi = ( A - \lambda^2 ) \bar \psi \,, \\
 & i \partial_t \varphi = ( A + \lambda^2 ) \varphi \,, \ \
i \partial_t \bar \varphi = ( A - \lambda^2 ) \bar \varphi \,. 
\end{split}
\end{equation}

The solutions $ \psi $ and $ \varphi $ have analytic extensions in 
$ \lambda $ to the upper half plane and $ \bar \psi $ and $ \bar \varphi $ 
to the
lower half plane. Same is true for $ a( \lambda ) $ and 
$ \overline{ a ( \lambda ) } $ respectively. Except in very special 
cases (such as our potential $ \alpha \, \sech x $) $ b ( \lambda ) $
does not have an analytic extension off the real axis. The 
reflection coefficient is defined as 
\[  r ( \lambda ) = \frac{ b ( \lambda ) } { a ( \lambda )}\,.\]

We assume that $ a ( \lambda ) $ has at most one zero and that it can 
only lie in $ \Im \lambda > 0 $, $ \Re \lambda = 0 $. That zero, 
$ \lambda_0 $, corresponds to an $ L^2 $ eigenfuction of $ L $, and 
at $ \lambda = \lambda_0 $, the two solutions 
are proportional:
\begin{equation}
\label{eq:defga}
 \varphi ( x , \lambda_0 ) = \gamma_0 \psi ( x , \lambda_0 ) \,. 
\end{equation}
The scattering data is given by the triple 
\begin{equation}
\label{eq:datab}
u ( t , x ) \longmapsto \{ r ( \lambda , t) \,, \lambda_0 \,, 
\gamma_0 ( t ) \} \,.
\end{equation}

The evolution of the scattering data 
is easily obtained from \eqref{eq:moreg}:
\begin{equation}
\label{eq:evol1b}
r ( \lambda , t ) = e^{ 2 i t \lambda^2} r ( \lambda , 0 ) \,, \ \ 
\gamma_0 ( t ) =  e^{ 2 i t \lambda^2} \gamma_0 ( 0 ) \,.
\end{equation}
In fact, we can 
use \eqref{eq:relb} \eqref{eq:moreg} to see
\[ \begin{split}  ( A - \lambda^2 ) \varphi &  =  
i\partial_t \varphi   = i\partial_t a  
 \bar \psi  + a    i \partial_t \bar \psi   + 
i \partial_t b   \psi   + b  
i \partial_t \psi \\ &  = 
i \partial_t a   \bar \psi+ ( i \partial_t 
b  + 2 \lambda^2 b )  \psi   + ( A - \lambda^2) \varphi   
  \,.
\end{split} \]
Independence of $ \psi $ and $ \bar \psi $ shows that, remarkably,
\[ \partial_t a ( \lambda, t ) = 0 \,, \ \ \partial_t b ( \lambda , t ) 
= 2 i \lambda^2 b ( \lambda , t) \,, \]
which gives the first part of \eqref{eq:evol1b}.
From \eqref{eq:defga} we see that,
\[   \gamma_0 ( A - \lambda^2 ) \psi ( \lambda_0 ) =  
i \partial_t \varphi ( \lambda_0 ) 
=  i \partial_t \gamma_0 
\psi ( \lambda_0 ) + \gamma_0 ( A + \lambda^2) \psi ( \lambda_0 )  \,, \]
so that, 
\[  \partial_t \gamma_0 ( t )  = 2 i \lambda^2 \gamma_0 ( t ) \,. \]
That gives \eqref{eq:evol1b}.
The justification of this formal calculation depends on 
$ u ( t , \bullet ) \in {\mathcal S} ( \RR ) $ and we refer to, for
instance, \cite[Section I.7]{FT} for a full proof.

\subsection{The Riemann Hilbert problem}
It is now universally acknowledged that the best way to obtain 
long time asymptotics for the inverse of
\eqref{eq:datab} and \eqref{eq:evol1b}
is by solving 
a Riemann-Hilbert problem  \cite{DIZ}, \cite[Chapter II]{FT}.
To recall this method let us consider the following 
matrix valued function of $ \lambda \in \CC \setminus \RR $, depending 
parametrically on $ x \in \RR $:
\begin{equation}
\label{eq:RH0}  \Psi ( \lambda, x ) \defeq \left\{ \begin{array}{ll}
{ [} a (  \lambda )^{-1}  \varphi ( x , \lambda ) 
e^{  i \lambda x } , 
 \psi ( x, \lambda  ) e^{- i \lambda x } ], & \Im \lambda > 0 \\ 
\ & \ \\
{ [ } \bar \psi ( x , \lambda ) e^{  i \lambda x } , 
\bar a( \bar \lambda)^{-1} \bar \varphi ( x , \lambda ) 
e^{ - i \lambda x } {]},
 & \Im \lambda <  0 \end{array}
\right. \end{equation}
The properties of $ \psi $ and $ \varphi $ (see for instance 
\cite[Section I.5]{FT}) imply that 
\begin{equation}
\label{eq:RH1}
 \Psi ( \lambda  ) = I + {\mathcal O} ( | \lambda|^{-1}  ) \,, \ \
|\lambda | \rightarrow \infty \,,  \end{equation}
where the decay rate may depend on $ x$, uniformly in compact sets.
From \eqref{eq:relb} we see that 
\[ \begin{split}
&   \Psi ( \lambda + i 0 , x )  = 
[ \bar \psi  ( x , \lambda ) 
e^{  i \lambda x } , 
 \psi ( x, \lambda  ) e^{- i \lambda x } ] \left[ \begin{array}{ll} 
\ \ \ \ \ \ \ 1 & 0 \\
 r ( \lambda ) e^{ 2 i ( \lambda x + \lambda^2 t ) } & 1 \end{array}
\right] \,, \\ 
&   \Psi ( \lambda - i 0 , x )  = 
[ \bar \psi  ( x , \lambda ) 
e^{  i \lambda x } , 
 \psi ( x, \lambda  ) e^{- i \lambda x } ] \left[ \begin{array}{ll} 
1 & -  { \bar r ( \lambda ) } e^{ - 2 i ( \lambda x + \lambda^2 t ) }  \\
0 & \ \ \ \ \ \ \ \ 1 \end{array}
\right] \,, \ \ \lambda \in \RR \,. \end{split} \]
Hence, the boundary values of $ \Psi ( \lambda ) $
satisfy
\begin{gather}
\label{eq:RH2}  
\begin{gathered}
\Psi ( \lambda + i 0 ) = \Psi ( \lambda - i 0 ) V_{x,t} ( \lambda ) 
\,, \ \
\lambda \in \RR \,, \\
V_{x,t} \defeq \left[ \begin{array}{ll}
\ 1 + | r ( \lambda )|^2   & { \bar r ( \lambda ) } e^{ - 2 i ( \lambda x + \lambda^2 t ) } \\ 
 r ( \lambda ) e^{ 2 i ( \lambda x + \lambda^2 t ) } & 
\ \ \ \ \ \ \ 1\end{array} \right] \,. 
\end{gathered}
\end{gather}
If $ a ( \lambda ) $ has no zeros in $ \Im \lambda \geq 0 $ then
the Riemann-Hilbert problem is to construct $ \Psi ( \lambda ) $ 
satisfying \eqref{eq:RH1} and \eqref{eq:RH2}. Liouville's theorem
readily shows that it is unique. If $ a ( \lambda ) $ has a zero
in $ \Im \lambda  > 0 $, and in our presentation we allow at most one, 
$ \lambda_0 $, we have to consider a Riemann-Hilbert problem in which
$ \Psi ( \lambda ) $ is allowed to have singularities at $ \lambda_0 $
and $ \bar \lambda_0 $. The structure of that singularity can be
seen in \eqref{eq:RH0}:
\begin{equation}
\label{eq:resb}
\begin{split}
&  \Res_{ \lambda = \lambda_0 } \Psi = \Psi ( \lambda_0 ) \left[ 
\begin{array}{ll} \ \ \ \  0 & 0\\
e^{ 2 i \lambda_0 x } \gamma_0' & 0 \end{array} \right] \,, \\
& \Res_{ \lambda = \bar \lambda_0 } \Psi = \Psi ( \bar \lambda_0 ) \left[ 
\begin{array}{ll} 0 & 
e^{ - 2 i \bar{\lambda_0} x } \bar \gamma_0'  \\
0 & \ \ \ \ 0 \end{array} \right] \,, \ \
 \gamma_0' \defeq \frac{\gamma_0 }{ a' ( \lambda_0 ) }\,.
\end{split}
\end{equation}
Since $ a ( \lambda ) $ can be reconstructed from $ r ( \lambda ) $
and $ \lambda_0 $ (see for instance \cite[Chapter I, (6.23)]{FT}; 
in our case it will be explicit) the Riemann-Hilbert problem in the 
case of one singularity is to find $ \Psi $ which in addition to 
\eqref{eq:RH1} and \eqref{eq:RH2} satisfies \eqref{eq:resb}.

A standard way to read off $ u ( t , x) $ from $ \Psi ( \lambda , x) $
follows from high frequency asymptotics of $ \psi ( x, \lambda )$ 
(see for instance \cite[(18)]{ZS72}):
\begin{equation}
\label{eq:asyb} \psi ( x , \lambda ) e^{ - i \lambda x } \simeq 
\left[ \begin{array}{l} 1 \\ 0 \end{array} \right] +
\frac{1}{ 2 i \lambda } \left[ \begin{array}{l} \ \ \ \ \  u (t,  x )  \\ 
\int_x^\infty | u ( t, y ) |^2 dy  \end{array} \right] + 
{\mathcal O} \left( \frac1 {| \lambda|^2 } \right) \,,
\end{equation}
so that 
\begin{equation}
\label{eq:recon}
  u ( x , t ) = \lim_{\lambda \rightarrow \infty } 2 i \lambda \Psi_{12} (
\lambda , x ) \,.
\end{equation}

We conclude this brief review by describing a reduction of the 
problem with prescribed singularities \eqref{eq:resb} to a problem
with $ \Psi $ analytic in $ \CC \setminus \RR $. To do that we
follow \cite[Section II.2]{FT} by considering a 
reformulation of the Riemann-Hilbert problem:
\begin{gather*}
G_{ \pm } ( \lambda )^{\mp 1 } \defeq \Psi ( \lambda) \left[
\begin{array}{ll} 1 & \ \ \ 0 \\ 0 & a_{\pm}( \lambda ) 
^{\mp 1 } \end{array} \right] 
\,,\ \ a_+ ( \lambda ) \defeq  a ( \lambda) \,, \ \ 
a_- ( \lambda ) \defeq \bar a ( \bar \lambda ) \,,  \ \  \pm \Im \lambda > 0 \,, \\
G_+ ( \lambda + i 0 ) G_-( \lambda -i0 ) = G ( \lambda ) \,, \ \ 
G ( \lambda ) \defeq \left[ \begin{array}{ll} 
\ \ \ 1 & - \bar b ( \lambda ) \\ - b ( \lambda ) & \ \ \ 1 \end{array}
\right] \,. 
\end{gather*}
The operators $ G_\pm ( \lambda ) $ are now analytic in $ \Im \lambda 
\pm > 0 $ 
(in fact, $ G_+^* ( \lambda ) = G_- (  \bar \lambda) $) 
but their ranks drop precisely at $ \lambda = \lambda_ 0 $
and $ \bar \lambda_0 $ respectively. The condition \eqref{eq:resb} 
becomes 
\begin{equation}
\label{eq:resb2} \Im G_+ ( \lambda_0 ) = {\rm span}_\CC \left[ \begin{array}{l}
\ \ \ \ 1 \\  \gamma_0 e^{ 2 i \lambda_0 x } \end{array} \right] \,, \ \
 \kerr G_- ( \bar \lambda_0 ) = {\rm span}_\CC \left[ \begin{array}{l}
- \bar \gamma_0  e^{ - 2 i \bar \lambda_0 x } 
\\ \ \ \ \ 1 \end{array}\right] \,.
\end{equation}
We now look for $ B ( \lambda ) $, analytic in 
$ \Im \lambda >  0 $, with $ B( \lambda )^{-1} $ analytic in 
$ \Im \lambda < 0 $, $ B ( \lambda ) = I + {\mathcal O} ( 1/|\lambda | ) $,
and such that 
\[ \widetilde G_+ ( \lambda ) \defeq G_+ ( \lambda ) B ( \lambda )^{-1} \,, \ \
\widetilde G_- ( \lambda ) = B ( \lambda ) G_+ ( \lambda ) \,, \ \ 
\pm \Im \lambda > 0 \,, \]
are nonsingular matrices. We note this requires $ B( \lambda) ^{-1} $
to have a pole at $ \lambda_0 $ and $ B( \lambda ) $, at $ \bar \lambda_0 $.
That is natural since we are adding to the ranks of $ G_\pm ( \lambda ) $. 
The condition \eqref{eq:resb2} mean that 
\begin{equation}
\label{eq:resb3} 
 \Im B ( \lambda_0 ) =  {\rm span}_\CC \widetilde G_+^{-1} ( \lambda_0 ) 
 \left[ \begin{array}{l}
\ \ \ \ 1 \\  \gamma_0 e^{ 2 i \lambda_0 x } \end{array} \right] \,, \ \
\kerr B ( \bar \lambda_0 ) =  {\rm span}_\CC \widetilde G_- ( \bar 
\lambda_0 )  \left[ \begin{array}{l}
- \bar \gamma_0  e^{ - 2 i \bar \lambda_0 x } 
\\ \ \ \ \ 1 \end{array}\right] \,. \end{equation}
This determines $ B ( \lambda) $ uniquely as a {\em Blaschke-Potapov
factor}:
\begin{gather}
\label{eq:bp}
\begin{gathered}
  B (\lambda ) = I + \frac{ \bar \lambda_0 - \lambda_0 }{ 
\lambda - \bar \lambda_0 } P \,, \ \ P^* = P \,, \ \ P^2 = P \,, \\
P = \frac{1}{  1 + |\beta|^2 } \begin{bmatrix} | \beta|^2 & \beta\\
 \bar \beta &  1 \end{bmatrix}\,, \ \
\beta ( x ) = \frac{  \widetilde G_-^{(11)} ( \bar \lambda_0 , x )  
\gamma_0 e^{ - 2 i \bar \lambda_0 
x } + \widetilde G_{-}^{(12)} ( \bar \lambda_0 , x ) } 
{ \widetilde G_-^{(21)} ( \bar \lambda_0 , x )  \gamma_0 e^{ - 2 i \bar \lambda_0 
x } + \widetilde G_{-}^{(22)} (\bar  \lambda_0 , x ) } \,,\end{gathered}
\end{gather}
see \cite[Chapter II, (2.17)-(2.27)]{FT}.

Hence to solve the Riemann-Hilbert problem \eqref{eq:RH2},\eqref{eq:resb}
we first solve the problem 
\begin{equation}
\label{eq:pstb} \widetilde \Psi ( \lambda + i 0 ) = 
\widetilde \Psi ( \lambda - i 0 ) V_{x,t} ( \lambda ) \,, \ \
\widetilde \Psi = I + {\mathcal O} ( 1/ | \lambda | ) \,, \ \ 
\text{ $ \widetilde \Psi $ analytic in $ \CC \setminus \RR $}\,. 
\end{equation}
We then define 
\[  \widetilde G_- ( \lambda ) = \widetilde \Psi ( \lambda ) 
 \left[
\begin{array}{ll} 1 & \ \ \ 0 \\ 0 & \bar a( \bar \lambda ) 
^{\mp 1 } \end{array} \right]\,, \ \
\widetilde G_+ ( \lambda ) = \widetilde G_- ( \bar \lambda ) ^* \,, \ \
\Im \lambda < 0 \,, \]
from which, using \eqref{eq:bp}, we construct $ B ( \lambda ) $.
Then 
\begin{equation}
\label{eq:resb4}
  \Psi ( \lambda ) = B ( \lambda )^{-1} \widetilde \Psi ( \lambda ) \,,
\end{equation}
and we can finally use \eqref{eq:recon} to obtain $ u ( t, x ) $.
In particular the long time behaviour of $ u ( t , x ) $ is determined by
the longtime behaviour of $ \widetilde \Psi ( \lambda ) $.

\subsection{Manakov ansatz}

The basic structure of the long time behaviour of $ \widetilde 
\Psi $ solving \eqref{eq:pstb} can be obtained from the 
{\em Manakov ansatz} for the solution of \eqref{eq:pstb} -- see
\cite[Section 2]{DIZ} and references given there. To describe it
we define
\begin{gather*}  m_\pm ( \lambda , x, t ) \defeq
 \frac{1}{ 2 \pi i } \int_\RR \frac{ r_\pm ( \zeta) \delta  ( \zeta+ i 0 )^{
\pm 1 } \delta ( \zeta -i0 )^{ \pm 1 } }{ \zeta - \lambda } 
e^{ \mp 2 (  i t \zeta^2 + x \zeta ) } d \zeta \,, \end{gather*}
where $ r_+ ( \zeta ) = r ( \zeta ) $, $ r_- ( \zeta ) = 
\bar r ( \zeta ) $, and 
\begin{gather*}
\delta ( \lambda, x, t ) 
\defeq \exp \left\{ \frac 1 {2 \pi i } \int_{-\infty}^{-x/2t}
\frac{ \log ( 1 + | r ( \zeta )|^2 ) } { \zeta - \lambda } d \zeta
\right\} \,, 
\end{gather*}
solves a scalar Riemann-Hilbert problem,  
\begin{gather*}
\delta ( \zeta + i 0 ) = \left\{ \begin{array}{ll} \delta ( \zeta 
- i 0 ) ( 1+ | r ( \zeta)|^2 ) \,, & \zeta < -x/2t \\
\ \ \ \delta( \zeta - i 0 ) \,, & \zeta > -x/2t \,, \end{array}
\right.  \end{gather*}
see \cite[Proposition 2.12]{DZ} for a detailed list of properties
of $ \delta ( z ) $ (stated in the defocusing case $ \log ( 1 - 
| r ( \zeta ) |^2 ) $). 
The Manakov ansatz is then given by 
\begin{equation}
\label{eq:manak}
\widehat \Psi ( \lambda , x,  t) \defeq 
\begin{bmatrix} \ \ \ 1 & m_+ ( \lambda , x , t ) \\
m_- ( \lambda , x , t) & \ \ \ 1 \end{bmatrix} 
\begin{bmatrix} \delta ( \lambda , x , t ) & \ \  1 \\ 
 \ \ 1 &  \delta ( \lambda , x , t )^{-1} \end{bmatrix} \,. 
\end{equation}
To see the properties of $ \widehat \Psi $ we use the following 
elementary lemma:
\begin{lem} 
\label{lb:1} 
Suppose that $ f \in C^\infty ( ( 0 , \RR ) ) $, $ f' \in L^1 ( [ 0 ,
\infty ) $, $ x^k f^{(l)} \in  L^\infty( [ 1 , \infty ) ) 
$, for all $ k$ and $ l$. 
Then, as $ \lambda \rightarrow \pm \infty $,
\[ \frac{1}{ 2 \pi i } \int_{\, 0}^\infty \frac{ f ( y) e^{ 
i \lambda  y^2 }} { y - x - i 0 } = 
\left\{ \begin{array}{ll} f ( x) e^{ i \lambda x^2 } + 
{\mathcal O} ( |\lambda|^{-\frac12} |x|^{-1} ) & 
\ \text{ $x > 0 $ and $ \lambda > 0 $} \\
\ & \ \\
{\mathcal O} ( |\lambda|^{-\frac12} |x|^{-1} ) & 
\ \text{ otherwise.} \end{array} \right. \]
\end{lem}
Using this lemma one checks that 
\[ \widehat \Psi ( \lambda - i 0, x , t )^{-1}\widehat \Psi ( \lambda + i 0 , x , t) 
= V_{x, t} ( \lambda ) + {\mathcal O} \left( \frac{1}{ \sqrt t 
| \lambda + x/2t | } \right) \,, \]
from which it follows as in \cite[Section 2]{DIZ} that 
\[ \widetilde \Psi ( \lambda ) = ( I + {\mathcal O} ( 1/\sqrt t) )
\widehat \Psi ( \lambda ) \,.\]
If we defined $ \tilde u ( t , x) $ by putting $ \widetilde \Psi  $ in
\eqref{eq:recon} this shows that $ \tilde u ( t ,x ) = 
{\mathcal O} ( 1 / \sqrt t ) $, and in fact a more precise 
statement can be obtained by using the second component in 
\eqref{eq:asyb}. 

We can now use \eqref{eq:bp} to obtain an approximation to 
$ B ( \lambda ) $ as $ t \rightarrow \infty $:
\[  B( \lambda ) = \widehat B ( \lambda ) + {\mathcal O} ( t^{-1/2} 
|\lambda|^{-1} ) \,, \ \ 
\widehat B ( \lambda ) = I -  \frac{ \bar \lambda_0 - \lambda_0 }{ 
\lambda - \bar \lambda_0 } \widehat P \,, \]
where $ \widehat P $ is as in \eqref{eq:bp} with $ \beta $ replaced by 
\[ \hat \beta = \gamma_0 \exp ( - 2 \bar \lambda_0 x ) 
\exp  \left\{ \frac 1 { \pi i } \int_{-\infty}^{\, 0}
\frac{ \log ( 1 + | r ( \zeta )|^2 ) } { \zeta - \lambda_0 } d \zeta
\right\} \,.\]
Hence to obtain a long time approximation for a solution we 
apply the procedure of \cite[Section II.5]{FT} to $ B_0 ( \lambda )$
since that corresponds to using \eqref{eq:recon} with 
\[ \Psi ( \lambda ) = B ( \lambda )^{-1} \widetilde\Psi ( \lambda ) 
= B ( \lambda)^{-1} ( I + {\mathcal O} ( t^{-1/2} ) ) \widehat \Psi 
( \lambda ) = B_0 ( \lambda )^{-1} ( I + {\mathcal O}(t^{-1/2}|\lambda|^{-1}
 ) ) \,.\]
This gives 
\begin{lem}
\label{lb:2}
Suppose that $ u ( 0 , \bullet ) \in {\mathcal S} ( \RR ) $ and that
the scattering data \eqref{eq:datab} for $ u ( 0 , \bullet ) $ is 
given by $ r ( \lambda ) $, $ \gamma_0 $ and $ \lambda_0  $, 
$ \Im \lambda > 0 $. Then 
\begin{gather}
\label{eq:lb2}
\begin{gathered}
u ( x, t ) = e^{ i \varphi_0 } \nlso \left( e^{ 2 i \Re \lambda_0 \, \bullet } 
2 \Im \lambda_0 \, \sech ( 2 \Im \lambda_0 \,  ( \bullet - x_0) ) \right) 
+ {\mathcal O}_{L^\infty}  ( t^{-1/2} ) \,, \ \ 
\end{gathered}
\end{gather}
where 
\[ x_0 = \frac{1}{ 2 \Im \lambda_0 } \left( \log | \gamma_0 | - 
\log | a' ( \lambda_0 ) | - \log (2 \Im \lambda_0 ) 
- \frac{\Im \lambda_0 }{ \pi } \int_{-\infty}^{\, 0 }
\frac{ \log ( 1 + | r ( \zeta )|^2 ) } { ( \zeta - \Re \lambda_0) 
^2 + \lambda_0^2 } d \zeta \right) \,, \]
and 
\[ \varphi_0 = \arg \gamma_0 - \arg a' ( \lambda_0 ) + \frac{1}{ \pi} 
 \int_{-\infty}^{\, 0 }
\frac{ \log ( 1 + | r ( \zeta )|^2 ) } { ( \zeta - \Re \lambda_0) 
^2 + \lambda_0^2 } ( \zeta - \Re \lambda_0 ) d \zeta  \,. 
\]
\end{lem}
We state this important result as a lemma to stress the fact
that a better error estimates seem available if more advanced methods
\cite{DZ},\cite{DIZ} are used.

\subsection{Scattering matrix}
We now apply Lemma \ref{lb:2} to obtain \eqref{eq:apb}. For that
we need to find the scattering data 
\eqref{eq:datab} for $ u ( 0 ,x ) = \alpha \sech x $. That is done
by a well known computation \cite{M81},\cite[Sect.3.4]{Maib} 
which reappears in many scattering theories, from the free
$S$-matrix 
in automorphic scattering, to Eckhardt barriers in quantum chemistry. 
We quote the results:
\[ a ( \lambda ) = \frac{ \Gamma ( \frac12 - i \lambda ) }
{\Gamma ( \frac 12 + \alpha - i \lambda ) \Gamma ( \frac12 - 
\alpha - i \lambda ) } \,, \ \ 
 b ( \lambda ) = i \frac {\sin \pi \alpha }{ \cosh \pi \lambda } \,, \ \ 
r ( \lambda ) = \frac{ b ( \lambda ) } { a ( \lambda ) } \,. \]
We note that in this special case $ b $ and $ r $ are meromorphic
in $ \CC $ (with infinitely many ``nonphysical'' poles).
Also, 
\[ \lambda_0 = 2 \alpha -1 \ \ \text{ if $ 1/2 < \alpha < 1 \,, \ $ \ and } \
\ \gamma_0 = b ( \lambda_0 ) = i \,.\]
We need to compute $ x_0 $ and $ \varphi_0 $. In general when $ u ( 0 ,x ) $
is even then $ x_0 = 0 $ by symmetry considerations. Here we see 
it by using \cite[Chapter II, (2.6)]{FT} which shows that
\[ \begin{split} \log| a' ( \lambda_0 )|  & 
=  \log( 2 \Im \lambda_0 ) + 
 \frac{ \Im \lambda_0}{ 2 \pi  }
 \int_{-\infty}^{\infty }
\frac{ \log ( 1 + | r ( \zeta )|^2 ) } { \zeta^2 + \Im \lambda_0^2 } 
d \zeta   \\
& = \frac{\Im \lambda_0}{\pi}  
 \int_{-\infty}^{0}
\frac{ \log ( 1 + | r ( \zeta )|^2 ) } { \zeta^2 + \Im \lambda_0^2 } 
d \zeta \,. \end{split} \]
Thus the formula in Lemma \ref{lb:2} results in $ x_0 = 0 $. 
To compute $ \varphi_0 $ we need to find the following integral
\[ 
\int_0^\infty \log \left( 1 + \frac{ \sin^2 \pi \alpha }{ \cosh^2 
 \pi \zeta  } \right) \frac{ \zeta}{ \zeta^2 + ( 2 \alpha -1)^2 } 
d \zeta \,, \ \ 1/2 < \alpha < 1 \,. \]

\end{document}